%% file: main.tex
\documentclass[12pt]{amsart}
\input head

\newcommand{\wl}{w_\lambda}

\author{Grigory Mikhalkin and Stepan Orevkov}
\title{Real algebraic knots and links of small degree}
\address{Universit\'e de Gen\`eve,  Math\'ematiques, Battelle Villa, 1227 Carouge, Suisse}
\address{Steklov Mathematical Institute,
ul. Gubkina 8, 119991, Moscow, Russia
\newline
IMT, Universit\'e Paul Sabatier, 118 route de Narbonne,
31062, Toulouse, France
\newline
National Research University Higher School of Economics,
Vavilova 7, Moscow, 117312 Russia}
\begin{document}
\hfill {\em To the memory of Sergei Duzhin}
\begin{abstract}
The paper gives a topological as well as rigid isotopy classification of smooth irreducible
algebraic curves in the real projective 3-space for the case when the degree 
of the curve is at most six and
its genus is at most one.
%5 and the genus is 1
%as well as for the case when their degree is 6 and their genus  is 0.
\end{abstract}
\maketitle

\footnote{Research is supported in part by the grants 141329, 159240, 159581 and the NCCR SwissMAP
project of the Swiss National Science Foundation (G.M.), and by RSF grant,
project 14-21-00053 dated 11.08.14. (S.O.).}

%\section{Introduction}
%We classify smooth isotopy types of knots and links formed by
%real algebraic curves in  $\rp^3$ of degree $d\le 6$. This problem is straightforward for $d\le 4$ and was solved
%by Johan Bjorklund \cite{Bj} in the case of $d=5$ for rational algebraic curves.

The subject of this paper is the problem of topological classification of smooth algebraic curves in $\rp^3$
when their genus and degree are fixed. The $\rp^2$ counterpart of this problem had
%(topological classification of smooth algebraic curves in $\rp^2$)
originated in the celebrated work \cite{Harnack} of Harnack,
was popularized by Hilbert in his famous list of problems \cite{Hilbert} and
consequently it was well-studied over the last century.
In the same time the corresponding topological classification in $\rp^3$ remains relatively unstudied.
We note that such classification is straightforward in the case when the degree is $4$ or less
as all relevant curves are contained either in a plane or in a quadric surface.
The classification in the case of degree 5 and genus 0 was obtained
by Johan Bjorklund \cite{Bj}.
Our paper continues this work by providing the classification in the case of degree 5 and genus 1 as well as 
for degree 6 and genus $\le 1$.

\section{Planar knots}
Let $\R K\subset\rp^3$ be a smooth irreducible real algebraic curve of degree $d$ and genus $g$
consisting of $l$ connected components.
This means that $\R K\subset\rp^3$ is given by a system
of homogeneous real polynomial equations in four variables so that 
the locus $\C K\subset\cp^3$ of complex solutions of the same system of equations
is an irreducible complex curve of genus $g$ smoothly embedded to $\cp^3$    
and homologous to $d[\cp^1]\in H_2(\cp^3)=\Z$.
Recall that by the Harnack inequality we have $l\le g+1$.
We assume that $\R K$ is non-empty, i.e. that $l\ge 1$.

%Note that if $d$ is sufficiently small with respect to $g$ then 
%$\R K$ must be contained in a quadric, or even in the plane.
%Namely, we have the following lemma.
%The following lemma is classical.
\begin{prop}
\label{lRR}
We have $g\le\frac{(d-1)(d-2)}{2}$.
If $g=\frac{(d-1)(d-2)}{2}$ then
there exists a hyperplane $\R H\subset\rp^3$
such that $\R K\subset\R H$. 
%If $g=\frac{(d-1)(d-2)}{2}$ there exists a hyperplane
If $2d-g+\iota(2d,g)< 9$ then
%If $d<3+g-\max\{2g-2-d,0\}$ then
there exists a quadric surface $\R Q\subset\rp^3$
such that $\R K\subset\R Q$.

Here $\iota(2d,g)$ is the maximal possible irregularity (the rank
of the first cohomology group) of a line bundle of degree $2d$
over a surface of genus $g$. In particular,
$0\le\iota(2d,g)\le \max\{0,2g-1-2d\}$. 
\end{prop}
\begin{proof}
Note that if $\R K$ is contained in a plane then $g=\frac{(d-1)(d-2)}{2}$
by the adjunction formula.
If $\R K$ is not contained in a plane in $\rp^3$ then
we may find a linear projection $\lambda:\R K\to\rp^2$ 
such that $\lambda(K)\subset\rp^2$ is a reduced singular planar curve of degree $d$.
To get such $\lambda$ we may use a projection from a point contained in
the line tangent to a generic point of $\R K\subset\rp^3$.
Thus $g<\frac{(d-1)(d-2)}{2}$.

The vector space of homogeneous quadratic forms in $\rp^3$ is 10-dimensional.
The restriction of such form to $\R K$ gives a section of the line bundle
of degree $d$ over $\R K$ associated to the projective embedding taken twice.
We get a linear map between two vector spaces. 
By the Riemann-Roch formula the dimension of the target vector space 
is not greater than $1+2d-g+\iota(2d,g)$ so the hypothesis of the lemma
ensures that the kernel is nontrivial.
The inequality 
$\iota(d,g)\le \max\{0,2g-1-d\}$ follows from Serre's duality as
$2g-2-d$ is the degree of the inverse bundle twisted by the canonical class
of the curve. 
\end{proof}

\begin{coro}\label{dle6}
If $d\le 6$ and $g> 2d-9$ then 
there exists a quadric surface $\R Q\subset\rp^3$
such that $\R K\subset\R Q$.
\end{coro}
\begin{proof}
By Lemma \ref{lRR} it suffices to check that $2d-g< 9$ and
$2d-g+2g-1-2d< 9$. The first inequality follows from our hypothesis
while the second one translates to $g<10$. 
Since $d\le 6$ this inequality holds unless $\R K$ is a planar sextic
(of genus 10), but then $\R K$ is contained in a reducible quadric surface.
\end{proof}

\begin{defn}
We say that $\R K\subset\rp^3$ is
a {\em planar} link, if it is isotopic
to a smooth (not necessarily algebraic) link $L\subset\rp^3$ that is
contained in a hyperplane $\R H\subset\rp^3$.
\end{defn}

Planar links are trivial from the topological viewpoint.
Namely, we have the following straightforward statement. 
\begin{prop}
Any two planar links with the same number $l$ of components and
the same parity of degree $d$ are isotopic.
\end{prop}
\ignore{
We denote with $p^l_0$ the isotopy type of a planar link of even degree and with 
$p^{l-1}_1$ the isotopy type of a planar link of odd degree (so that we have one non-trivial
and $l-1$ trivial components.

\begin{proof}
Indeed, if $d$ is even then 
it is isotopic to the disjoint union of $l$ unknotted circles
as any pair of nested ovals of $L\subset\rp^2$ can be turned into
a non-nested pair by an isotopy in $\rp^3$.

Similarly, if $d$ is odd then $\R K$ is isotopic to a planar curve
in $\rp^2\subset\rp^3$ consisting of
of a pseudoline and $l-1$ non-nested ovals.
\end{proof}
}

Suppose now that $\R K\subset\R Q$ for
%a smooth link
%$L$ isotopic to $\R K$ and
a quadric surface $\R Q\subset\rp^3$.
There are four topological types for a quadric $\R Q$.
\begin{itemize}
\item The case when $\R Q$ is reducible (or non-reduced) and thus contains a plane.
Then $\R K$ is planar, so that $g=\frac{(d-1)(d-2)}2$.
%The topology of the link $\R K$ is thus very simple.
%If $d$ is even then it is the disjoint union
%of up to $g+1$ unknotted and unlinked circles.
%If $d$ is odd then it is isotopic to the disjoint union
%of the line $\rp^1$ and 
%of up to $g$ unknotted and unlinked circles.
\item The case when $\R Q$ is ellipsoid. Then $d$ is necessarily even
and $\R K$ is planar.
% is a collection of up to $g+1$ unknotted and unlinked circles
%as in the planar case.
\item The case when $\R Q$ is a singular quadric surface.
If $\R K$ is disjoint from the singular point of $\R Q$ then
$d$ is even and $\R K$ is planar
while $g=(\frac d2-1)^2$.
%topologically the situation is similar to the previous one:
%$\R K$ is contained in a smooth real surface isotopic to
%the standard sphere in $\rp^3$. It is a collection of
%up to $g+1$ unknotted and unlinked circles that we have already
%seen in the planar case.
A component of $\R K$ containing the singular point of $\R Q$
must be isotopic to the line in $\rp^3$. All other components
must bound disks in the cone $\R Q$, so once again
the link $\R K$ must be planar, but $d$ has to be odd.
%we have the planar configuration, but for an odd degree $d$.
%Note that if $\R K$ (rather than a smooth knot isotopic to it)
%sits on $\R Q$ and $\R Q$ is a conic surface we also have two cases.
%If $\R K$ is disjoint from the apex of $\R Q$ then $\R K$ is isotopic
%to a knot on an ellipsoid (as we can smooth $\R Q$ to ellipsoid) and
%thus have the bidegree $(\frac d2,\frac d2)$ there.
%In this case
%$\R K$ can be obtained from a smooth planar
%curve $\R \tilde K\subset \rp^2$ of degree $\frac d2$
%as the image under a regular degree 1 map
%$\R\hat Q\to\R Q$.
%Here the surface $\R \hat Q$ is obtained from $\rp^2$ 
%by a blowup in a point outside of $\R \tilde K$.
%We have $g=\frac{(\frac d2-1)(\frac d2-2)}2$ in this case. 

Recall that $\R Q$ can be obtained from the (toric) Hirzebruch surface $F_2$
by contracting the $(-2)$-sphere there.
If $\R K\subset\R Q$ passes through the singular point of $\R Q$ then
$d$ must be odd. 
Furthermore, the curve $\R K$ is the image of a smooth curve $\R \tilde K\subset \R F_2$
whose Newton polygon is the trapezoid with vertices $(0,0),(0,\frac{d-1}2),
(1,\frac{d-1}2),(d,0)$. 
%$\R \hat K\subset\R \hat Q$ under a regular degree 1 map
%$\R\hat Q\to\R Q$. Here $\R \hat K\subset\R \hat Q$ is obtained from
%a smooth planar curve $\R \tilde K\subset\rp^2$ of degree $\frac{d+1}2$.
We have $g=(\frac{(d-1)(d-3)}4)$ if $d$ is odd, just as in the case of a curve
of bidegree $(\frac{d+1}2,\frac{d-1}2)$ on a hyperboloid.

\item The case when $\R Q$ is a hyperboloid. 
This is the most interesting case. We consider it in more details in the following
section.
\end{itemize}

\section{Hyperboloidal links}
\begin{defn}
We say that $\R K\subset\rp^3$ is
a {\em hyperboloidal} link, if it is isotopic
to a smooth (not necessarily algebraic) link $L\subset\rp^3$ that is
contained in the hyperboloid $\R Q=\{(x:y:z:u)\in\rp^3\ |\ x^2+y^2=z^2+u^2\}\subset\rp^3$.
\end{defn}

Such $L$ consists of $k$ components that bound disks in $\R Q$
and $j$ homologically non-trivial components.
Note that all non-trivial components must be homologous
in $\R Q$ as otherwise they would intersect.
Recall that the hyperboloid $\R Q$ has a distinguished
basis (up to a sign and permutation) in its homology group $H_1(\R Q)=\Z\oplus\Z$ given
by the ruling $\R Q=\rp^1\times\rp^1$.
%Thus we may associate to $L\subset\R Q$ 

To any hyperboloidal link $L$ we prescribe three integer numbers: $jp,jq,k$.
Here $(p,q)$ is the only non-trivial homology class of a component of $L$ in
$H_1(\R Q)=H_1(\rp^2\times\rp^2)=\Z\oplus\Z$ (so that $p$ and $q$ are coprime,
$j$ is the number of components in this 
class and $k$ is the number of homologically trivial components (ovals) of $L$ in $\R Q$.
Choosing the orientation of the generating lines of $\R Q=\rp^1\times\rp^1$
as well as their order in the basis we may assume that $p\ge q\ge 0$.

If $jp=jq$ then the non-trivial components of $L$ on $\R Q$ bound disjoint disks in $\rp^3$.
These disks can be obtained from planar sections of $\R Q$ and thus $L$ is planar in this case.

Consider the case $jp=jq+1$, so that $j=1$ and $p=q+1$.
The only non-trivial component of $L$ intersects a curve $S\subset\R Q$
of homology class $(1,-1)$ in a single point. In its turn, $S$ can be cut on $\R Q$
by a plane in $\rp^3$. We may contract $S$ to a point so that $\R Q$ becomes
a quadratic cone and $L$ becomes a link on this cone passing through its apex.
Thus $L$ is planar (of odd degree).
We have proved the following statement.

\begin{prop}\label{plan-hyp}
If $L$ is a non-planar hyperboloidal link then $p>q$.
Furthermore, if $p=q+1$ then $j>1$.
\end{prop}

\begin{defn}
We denote with $h_{a,b}\ \sqcup \langle k
\rangle$, $a>b+1$, $k\ge 0$,
the isotopy type of a hyperboloidal link $L$ with $j=GCD(a,b)$
non-trivial components on $\R Q$.
Here $(\frac aj,\frac bj )\in H_1(\R Q)=H_1(\rp^1\times\rp^1)=\Z\oplus\Z$
is the homology class of a component of $L$ (with the appropriate orientations and order
of the basis elements) and $k$ is the number of homologically trivial components (ovals).
We use the abbreviated notation $h_{a,b}$ for
$h_{a,b}\ \sqcup \langle k\rangle$ when $k=0$.
\end{defn}

\begin{prop}
The hyperboloidal links of type
$h_{a,b}\ \sqcup \langle k\rangle$, $a>b+1$, are non-isotopic for different
values of $a,b,k$.
\end{prop}
\begin{proof}
Consider the universal covering $\pi:{\mathbb S}^3\to\rp^3$. Let $\pi_{\R Q}:\R \tilde Q\to\R Q$
be the restriction of this double covering to $\R Q$.
The covering $\pi_{\R Q}$ is given by the subgroup
$\{(\alpha,\beta)\in H_1(\R Q)\ |\ \alpha+\beta\equiv 0\pmod 2\}\subset H_1(\R Q)=\pi_1(\R Q)$,
thus the total space $\R \tilde Q$ is a torus.
Furthermore, this torus is a standard torus in ${\mathbb S}^3$ (the boundary of a tubular
neighborhood of an unknot) and the classes $(1,1)$ and $(1,-1)$ in $H_1(\R Q)$
are the images of its standard generators (bounding disks in ${\mathbb S}^3\setminus\R\tilde Q$). 

Thus $\pi^{-1}(L)$ is the union of a $(a+b,a-b)$-torus link in ${\mathbb S}^3$ with $2k$ unknotted unlinked 
circles. Since by our hypothesis we have $a-b>1$, all such links are non-isotopic.
\end{proof}

\def\S{\mathbb S}
\def\T{\mathbb T}
\def\RP{\mathbb{RP}}
\def\eps{\varepsilon}

\begin{rmk}%[Hyperboloidal links and projectivization of toric links]
We may describe the relation between the hyperboloidal links
and toric links as follows.
Let $\S^3=\{|z|^2+|w|^2=1\}$ be the unit sphere in $\CC^2$ and let
$\T$ be the torus $\{|z|=|w|=1\}\subset\S^3$. As usually, for $p,q\in\Z$,
we define the  {\it$(p,q)$-torus link} as $T(p,q)=\{z^p=w^q\}\cap\S^3\subset\T$.
It is clear that $T(p,q)\sim -T(-p,q)\sim T(q,p)$ and it is well known that
$T(p,q)$ is determined by $(p,q)$ up to isotopy under the condition $p\ge|q|$.
The number of components of $T(p,q)$ is equal to $\gcd(p,q)$.

In the case when $p\equiv q\mod 2$, the link $T(p,q)$ is invariant under the antipodal
involution $-1:\CC^2\to\CC^2$, $(z,w)\mapsto-(z,w)$. So, in this case we define the
{\it projective $(p,q)$-torus link} as the quotient $\bar T(p,q)=T(p,q)/(-1)$. It
sits in $\S^3/(-1)$ which we naturally identify with $\RP^3$.
It is clear that the isotopy type of $\bar T(p,q)$ is determined by $(p,q)$ up to the
above relations. Indeed, if two links in $\RP^3$ are isotopic, then their double covers
are isotopic as well.

As we have already seen, we have $h_{a,b}=\bar T(a+b,a-b)$.
\end{rmk}

To represent toric and projective toric links by diagrams, it is convenient to use the language
of braids. We define the closure of a braid in $\S^3$ in the usual way and we define
the {\it closure in $\RP^3$} of a braid as follows. A braid can be naturally identified with
a tangle in a (round) ball $B^3$ with all endpoints placed symmetrically on a great circle on $\partial B^3$.
So, the closure of the braid in $\RP^3$ is the image of the tangle under the identification of
antipodal points of $\partial B^3$.
In particular, the diagram of the closure of a braid in $\RP^3$ just coincides with the diagram of the
corresponding tangle.

Let $p$ and $q$ be positive of the same parity and $\eps=\pm1$. Then $T(p,\eps q)$ is the closure in $\S^3$
of the $p$-braid $(\alpha\beta)^q$ where
$\alpha=\sigma_1^\eps\sigma_3^\eps\dots$ and $\beta=\sigma_2^\eps\sigma_4^\eps\dots$.
Similarly, $\bar T(p,\eps q)$ is the closure in $\RP^3$ of the braid represented by the first half of the
word $(\alpha\beta)^q$, i.~e., the braid $(\alpha\beta)^{q/2}$ if $p$ and $q$ are even and
$(\alpha\beta)^{(q-1)/2}\alpha$ if $p$ and $q$ are odd.

Many of the knots appearing in the classification
results of this paper are hyperboloidal (as topological knots)
even if the corresponding spatial algebraic curves
are not necessarily contained in a quadric surface.
E.g. we have
$K_3=h_{4,1}=\bar T(3,5)$ in Figure \ref{d5g0},
while we have
$K_5=h_{3,1}=\bar T(2,4)$,
$K_8=h_{5,3}=\bar T(2,8)$,
$K_{11}=h_{7,5}=\bar T(2,12)$,
$K_{14}=h_{5,1}=\bar T(4,6)$ in Figure \ref{deg6c}.
Note that among these knots, only $h_{4,1}$ and $h_{5,1}$ are realizable by
rational algebraic curves of respective degree (5 and 6) sitting in a hyperboloid. 

\section{Viro's invariant}
\input{vironew.tex}

\section{Projection from a double point
and the resulting diagram}
\input dpoint.tex

\section{Knots and links of degree up to 5}
\input upto5.tex

%\section{Knots and links of degree 6}

\input deg6.tex

\input d6g1.tex

\end{document}

%% file: head.tex
\usepackage{amsmath}
\usepackage{amssymb}
\usepackage{amsthm}
\usepackage{graphicx}

\newcommand{\ignore}[1]{\relax}

\newcommand{\Id}{\operatorname{Id}}

\newcommand{\C}{}

\newcommand{\R}{\mathbb R}
\newcommand{\Z}{\mathbb Z}

\newcommand{\T}{\mathbb T}

\newcommand{\ind}{\operatorname{ind}}
\newcommand{\lk}{\operatorname{lk}}

\newcommand{\CC}{\mathbb C}
\newcommand{\p}{\mathbb P}

\newtheorem{lem}{Lemma}
\newtheorem{claim}{Claim}[section]

\newtheorem{theorem}[lem]{Theorem}
\newtheorem{corollary}[lem]{Corollary}
\newtheorem{thm}{Theorem}

\newtheorem{coro}[lem]{Corollary}
\newtheorem{prop}[lem]{Proposition}

\theoremstyle{definition}
\newtheorem{condition}[claim]{Condition}
\newtheorem{defn}[lem]{Definition}

\newtheorem{problem}[claim]{Problem}

\theoremstyle{remark}
\newtheorem{rmk}[lem]{Remark}
\newtheorem{rem}[lem]{Remark}

\renewcommand{\setminus}{\smallsetminus}

\newcommand{\conj}{\operatorname{conj}}

\newcommand{\dd}{\partial}

\newcommand{\cp}{
{\mathbb P}}
\newcommand{\rp}{{\mathbb R}{\mathbb P}}

\newcommand{\pp}{{\mathbb P}}

\renewcommand{\setminus}{\smallsetminus}

\ignore{

\newtheorem{condition}[theorem]{Condition}

\newenvironment{proof}[1][Proof]{\noindent\textbf{#1.} }{\ \rule{0.5em}{0.5em}}
}

%% file: vironew.tex
As it was discovered in \cite{Vi}, to any real algebraic link $K\subset\p^3$
one can associate an integer invariant (called in \cite{Vi}
{\em encomplexed writhe}) which we denote with $w$.
We briefly recall its definition.

The projection $$\pi_p: K\to\cp^2$$ from a point
$p\in\rp^3\setminus\R K$
maps the
%$\R K$ into a planar curve $\R C\subset\rp^2$.
%Furthermore, it maps the
{\em complexification} $\C K\subset\cp^3$ (i.e. the set of complex points of $\R K$)
into a planar complex curve $\C C\subset\cp^2$.
If $p$ is chosen generically then $\C C$ is {\em nodal},
i.e. all its singularities are nodes, i.e. transverse intersections
of pairs of local branches of $\pi_p(\C K)$.

The nodal curve $\C C$ is defined over $\R$, but the set $\R C\subset\rp^2$ of its
real points may be different from $\pi(\R K)$
as $K$ is the normalization of $C$.
Clearly we have
$\pi_p(\R K)\subset\R C$. Each point $q\in \R C\setminus\pi_p(\R K)$
must be singular and thus is a node. Thus $\pi_p^{-1}(q)$ consists
of a pair of conjugate points of $\C K\setminus \R K$.
We refer to such a node $q\in\R C$ as {\em elliptic},
in contrast to the nodes of $\pi_p(\R K)$ which are called
{\em hyperbolic}.

We denote that set of nodes of $\R C$ with $\Sigma\subset\R C$.
As in the classical knot theory, 
the algebraic curve $\R C$ may be considered as a {\em knot diagram} of
$\R K$. 
The hyperbolic nodes of $\R C$ correspond to conventional self-crossings
of classical knot diagrams.
But $\R C$ may also contain elliptic nodes, visible in $\rp^2$
as solitary points disjoint from the rest of the real curve.
We treat the elliptic nodes from $\Sigma$
as imaginary self-crossings points of the diagram $\R C$.

\begin{figure}[h]
\includegraphics[height=12mm]{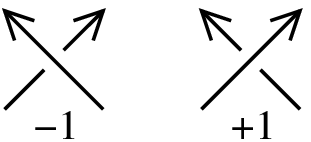}
\caption{Signs of diagram self-crossing points.\label{writhe}}
\end{figure}

Suppose that a hyperbolic node $q\in\Sigma$ corresponds to an intersection
of two branches from the {\em same} component of $\R K$.
Then, in full consistency with the knot theory conventions,
we define the sign $\sigma(q)=\pm 1$ according to Figure \ref{writhe}
by choosing an affine chart of $RP^3$ such that $p$ is at the infinity.
One easily checks that the definition of the sign does not depend on
a choice of the affine chart.
\ignore{
Namely, we choose an orientation of $\R K$, choose an oriented
straight interval chord in $\rp^3$ connecting the two points of $\pi_p^{-1}(q)$
in a way disjoint from $p$,
and order the two points of $\pi_p^{-1}(q)$ according to the
chosen orientation of the chord.
Then we compare the orientation of the 3-frame formed by
the oriented tangent vectors to $\R K$ (in the given order) and the oriented
tangent vector to the chord against the standard orientation of $\rp^3$.
If the two orientations agree we set $\sigma(q)=+1$ and 
otherwise we set $\sigma(q)=-1$.
Clearly $\sigma(q)$ is independent of the choices we made:
the orientation of the chord appear twice and the orientation of
the component of $\R K$ containing  $\pi_p^{-1}(q)$ also appears twice
in the definition of this sign.
}

The paper \cite{Vi} has also introduced a similar sign for
the elliptic nodes. Namely, for an elliptic node $q\in\Sigma$
we have $\pi^{-1}_p(q)\in\C K\setminus\R K$.
Note that $p$ and the two points of $\pi^{-1}_p(q)$ sit on the same line (defined over $\R$)
which we
denote with
$l_q\subset\p^3$.
%passing through $p$ and $\pi_p(q)$.
Note also that $l_q\setminus \R l_q$ is
the disjoint union
of two open hemispheres of the Riemann sphere
$l_q=\cp^1$.
%(canonically oriented by the complex orientation)
%that is the complexification of $l_q$.
As the two points of
$\pi^{-1}_p(q)$ are conjugate, they
belong to different hemispheres.
Choose a component $S_l$ of $l_q\setminus \R l_q$.
The hemisphere $S_l$
is canonically oriented as an open set of a complex curve.
Thus the choice of $S_l$ defines an 
orientation of $\R l_q$ (through $\R l_q=\dd S_l$),
as well as the point $u\in\pi^{-1}_p(q)\cap S_l$.
%contained in that component.
The orientation of $S_l$ induces
a local orientation of $\rp^2$ at $q$ with the help of
$d\pi_p|_{u}:T_{u} S_l\to T_q \rp^2$.
Together with the orientation of $\R l_q$ it defines an orientation
of $\rp^3$. The sign $\sigma(q)=\pm 1$ is defined by comparison
against the standard orientation of $\rp^3$.
As before, it does not
depend on the auxiliary choice of orientation of $\R l_q$ as
this choice enters the definition of $\sigma(q)$ twice.

For simplicity of notation we define $\sigma(q)=0$ if
$q\in\Sigma$ is a hyperbolic node corresponding to intersections
of different components. % of the normalization of $\R K$.

\begin{thm}[Viro \cite{Vi}]\label{virothm}
The sign $\sigma$ at hyperbolic and elliptic
nodes of $\R C$ is consistent in the following sense.
% that the nodal set $\Sigma_t$
% under the first Reidemeister move,
%see Figure \ref{reidem1}.
For any 1-parametric family $\R K_t\subset\rp^2$,
$0\le t\le 1$ of embedded smooth algebraic curves
of the same degree %and the same geometric genus 
with $p\notin\R K_t$ the union of the corresponding sets $\Sigma_t$ of real nodes
of $\pi_p(\C K_t)$ enhanced with the signs $\sigma$
defines an integer homology 1-chain with the boundary $[\Sigma_1]-[\Sigma_0]$. 
In other words, when $t$ varies only pairs of nodes of opposite signs may annihilate
while remaining nodes keep their sign invariant.
\end{thm}

A 1-parametric family of smooth algebraic curves of the same degree
is also called {\em rigid isotopy}.
\ignore{
If our rigid isotopy $\R K_t\subset\rp^3$ is a smooth and
generic path in the space of embedded curves of given degree
then the diagrams $\R C_t\subset\rp^2$ are nodal for all but
finitely many values of $t\in (0,1)$.
Clearly all nodes together with their signs are preserved
under rigid isotopies such that all diagrams $\R C_t$ are nodal.

The paper \cite{Vi} identified the 5 types
of generic non-nodal diagrams
$\R C_t\subset\rp^2$.
Namely, $\R C_t$ can be cuspidal, then the passage between $\R C_{t-\epsilon}$
and $\R C_{t+\epsilon}$ corresponds to the first Reidemeister
move, where a hyperbolic node on one side becomes an elliptic node
of the same sign on the other side.
}
\begin{figure}[h]
\includegraphics[height=12mm]{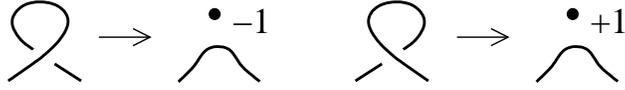}
\caption{Interchange between elliptic and hyperbolic nodes of the same sign.\label{reidem1}}
\end{figure}
\ignore{
Also under a generic rigid isotopy $\R C_t$ can develop a tacnode.
Then a pair of real nodes of different signs can disappear into 
a pair of complex conjugate nodes in
$\CC\p^2\setminus\rp^2$.
There are two types here which correspond to annihilation
of a pair of elliptic, and a pair of hyperbolic nodes.
This correspond to two real versions of the second Reidemeister move.

Finally, under a generic rigid isotopy $\R C_t$ can develop
an ordinary triple point. Then $\R C_{t-\epsilon}$
and $\R C_{t+\epsilon}$ have the same number and type of nodes.
There are two types here as well as an ordinary triple point
defined over $\R$ may consist of three real branches, or
of one real and a pair of complex conjugate branches.
These types correspond to the third Reidemeister move.
}
\begin{coro}[Viro \cite{Vi}]
The number
\begin{equation}\label{w}
w(K)=\sum\limits_{q\in\Sigma}\sigma(q)
\end{equation}
is invariant under the rigid isotopy of the real algebraic curve $\R K$ as
well as the choice of the point $p\in\rp^3\setminus\R K$ determining
the diagram $\R C\subset\pi_p(\C K)$.
\end{coro}
\ignore{
Invariance under rigid isotopies follow from Theorem \ref{virothm}.
A deformation $p_t\in\rp^3$, $t\in[0,1]$
of the projective point corresponds to a rigid
isotopy of curves after applying projective linear 
automorphisms of $\rp^3$ sending $p_t$ to $p_0$.
}
%If $K$ is of type I then it is also convenient
%to consider the number
%$$\wl(K)=w(K)+\sum\limits_{M,L 

%% file: dpoint.tex
Consider the divisor $D\subset\C K$ obtained
by intersecting $K$
with a generic plane in $\p^3$.
All effective divisors linear equivalent to $D$ form
a linear projective space $|D|\approx \pp^r$. %of dimension $r$.
The Riemann-Roch theorem ensures that $r\ge d-g$.

\begin{figure}[h]
\includegraphics[height=12mm]{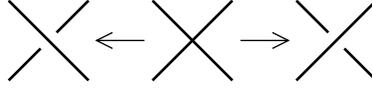}
\caption{Two resolutions (perturbations) of a self-crossing point
of the curve.
\label{perturb}}
\end{figure}

\begin{lem}\label{rkprime}
If $r\ge 4$ then there exists a continuous deformation
$$f_t:K\subset\p^3,\ 0\le t\le 1,$$
in the class of real algebraic curves of degree $d$
(keeping the source $K$ unchanged as an abstract real curve)
with the following properties.
\begin{itemize}
\item For $0\le t< 1$ the map $f_t$ is 
an embedding, and $f_t(K)$ is
a smooth real algebraic curve of degree $d$.
%whose complexification is smooth.
\item $f_0=\Id$.
\item The map $f_1:K\to\p^3$ is an immersion with
a single self-crossing point $p$,
with $D^-=f_1^{-1}(p)$ consisting of two points.
%$x,y\in \C K$.
We have $D^-=\conj D^-$.
%The pair $\{x,y\}$ is real,
%i.e. either $x,y\in\R K$ or $\conj(x)=y$.
%The restriction $f_1|_{\C K\setminus\{x,y\}}$
%is an embedding.
The two branches of $f_1$ at $p$ are
not tangent to each other. 
\end{itemize}

Vice versa, if $r\ge 4$ then any immersed algebraic curve
with a single self-crossing point
of two (real or imaginary) branches with distinct tangent directions
can be perturbed to an embedded real algebraic
knot obtained by perturbing the two branches
in any of the two possible directions, i.e. so that each
sign of the resulting crossing point on the diagram may appear
%the two intersection points 
%of the cycle \eqref{cycle} and the Schubert cycle of lines passing
%through a point of $\rp^3$ sufficiently far from $p$
(see Figure \ref{perturb} for the case when the two branches are real). 
\end{lem}

In other words,
$f_t$ is an isotopy of $K$ to
a curve with a single self-crossing point.
Conversely such a point can be resolved in
any direction.
\begin{proof}
The curve $K\subset\p^3$ is given by a choice of
a 3-dimensional projective subspace in $|D|$.
Consider $\tilde K\subset\rp^4$ corresponding
to enlarging this subspace to a 4-dimensional linear
subspace in $|D|$, so that we have
$K=\pi_q(\tilde K)$ for the projection $\pi_q$ from
$q\in\rp^4\setminus\R \tilde K$. 
%$$
%\pi_q:\rp^4\setminus\{q\}\to\R P^3
%$$
%from a point $q\in\rp^4\setminus\R \tilde K$.

Consider the subspace $\C\Xi\subset\cp^4$
formed by all projective chords of $\C \tilde K$, i.e. lines
connecting two points of $\tilde K$ (or the tangent
line in the case when the two points coincide).
The subspace $\C\Xi$ is a real algebraic hypersurface.
Note that projections of $\tilde K$ from generic
points of $\R\Xi$ are immersed curves in $\p^3$
as the tangent lines to $\C\tilde K$ form a 2-dimensional
stratum in $\C\Xi$.

Suppose that $\pi_q(\C\tilde K)\subset\p^3$ has  
a self-crossing point $p\in\p^3$ with
tangent branches.
Then there exists a plane in $\cp^4$ tangent
to $\C\tilde K$ at $\pi_q^{-1}(p)$ . 
However, for each point of $\C \tilde K$
there might be only a finite number of other points 
of $\C \tilde K$ sharing a tangent plane.

Thus, an interval $q_t$, $0\le t\le 1$,
in $\rp^4$ connecting $q$ to
a generic point of $\R\Xi$ (chosen in the boundary
of the component of $\rp^4\setminus\R\Xi$ containing $q$)
produces the requited deformation
$f_t(K)=\pi_{q_t}(\tilde K)$.

Conversely, moving a generic point of
the hypersurface $\R\Xi$
to any of the two sides we get a perturbation
of the map with a self-crossing into an embedding.
\end{proof}

Consider the curve $K'=f_1(K)$ corresponding to 
a generic point $q\in \R\Xi$ from the 
%obtained from $\R K$ with the help of
%from 
proof of Lemma \ref{rkprime}.
This curve is singular and the map
$f_1:K\to K'$ can be viewed as the normalization.
Consider the planar curve
$$C=\tilde\pi_p(K)\subset \p^2,$$
where $\tilde\pi_p$%:\R \tilde K'\to\rp^2$$
%be  the map from the normalization
%$\R \tilde K'$ of $\R K'$
is the lifting
of the projection $\pi_p|_{K'\setminus\{p\}}$
to $K$.
%obtained through
%extension of the projection
%$\pi_p:\rp^3\setminus\{p\}$
%from the double point $p\in\R K'$
%restricted to $\R K'\setminus\{p\}$.
%Denote $\R C=\tilde\pi_p(\R K')$.
%(note that $\pi_p$ is multivalued at $p$ 
%as usual $\pi_p(\R K')$ denotes the entire image).

Recall that a planar curve is called {\em nodal}
if all of its singular points are simple nodes. 
\begin{prop}\label{pirkprime}
The curve $C\subset\p^2$ is a nodal
curve of (geometric) genus $g$ and degree $d-2$
in $\p^2$. 
Furthermore, $\tilde\pi_p(D^-)$ is disjoint
from the nodes of $C$.
\end{prop}
\begin{proof}
The geometric genus of $C$ coincides with that of
its normalization $K$.
%as the two curves are parameterized by the same 
%Riemann surface.
The degrees of $C$ in $\p^2$ and of $K'$
in $\p^3$ differ by two since the inverse image of
a line in $\cp^2$ is a plane that intersect $K'$
at $p$ with the local multiplicity 2.

The curve $C$ is obtained by a projection of
$\tilde K\subset\rp^4$ from a line $l\subset\rp^4$\
passing through two points of $\tilde K$.
% $x,y\in\R\tilde K$. 
The condition that $C$ is an immersion is equivalent
to the condition that no plane in $\rp^4$ containing
$l$ can be tangent to $\tilde K$.

The space of planes tangent to $\tilde K$
is parameterized by $\R\tilde K$ itself and 
thus 1-dimensional. Each such
plane intersects $\tilde K$ in finitely many points.
As the number of pairs of points in $\tilde K$ is 
2-dimensional, $C$ is an immersion for a generic
choice of $q\in\R\Xi$.

Moving the chord $l$ of $\tilde K$
produces a 2-parametric family of deformations of $C$.
Away from a small neighborhood of $p$
it infinitesimally corresponds to 
moving the projection point $p\in\rp^3$
in the 2-plane tangent to the
branches of $\R K'$ at $p$. 
%Thus a generic choice of $p$ ensures
%a generically immersed curve $\R C$ with $x'',y''$
%disjoint from the double points.
%Outside a small neighborhood of $p$ this deformation
%is approximated by changing the projection at
\end{proof}

\ignore{
Recall that an effective divisor $\C Z$ on $\C C$
can be considered as
%a finite multisubset of $\C C$, i.e.
a finite subset of $\C C$ where the points are
prescribed multiplicities (natural numbers).
Somewhat loosely we write nevertheless $\C Z\subset\C C$.
The degree of $\C Z$ is the sum of multiplicities
of all points. We say that $\C Z$ is real if it is invariant
with respect to complex conjugation. In this case we
denote the fixed (multiset) locus of $\conj$ with $\R Z\subset\R C$.
Choose a plane $\R H\subset\rp^3$ disjoint 
from a fixed point $p\in\rp^3$.
}
%Let $D$ be a divisor on $K$. The following
%statement is straightforward.
%Recall that from the viewpoint of linear systems
%a real curve $C\subset\p^2$ normalized by $K$
%is a ..
From a real curve $C\subset\p^2$ normalized by
$\nu:K\to C$
we may reconstruct a spatial curve $K'\subset\p^3$
such that $C=\tilde\pi_p(K')$, $p\in K'$
with the help of a (not necessarily effective)
divisor $D$ on $K$ in the equivalence class
of a line section of $C$.
Namely, we have the following straightforward
statement. Let $D=D^+-D^-$ be a real
(i.e. invariant with respect to $\conj$)
divisor on $K$,
where $D^+$ and $D^-$ are disjoint effective divisors.
\begin{prop}
\label{reconstruct-rkprime}
There exists a curve $K'\subset \p^3$
normalized by $K$ (here we denote the normalization
$K\to K'$ with $f_1$ to make notations consistent with
Lemma \ref{rkprime}),
a point $p\in K'$ such that $\tilde\pi_p(K)=C$,
$f_1^{-1}(p)=D^-$, and a plane $H\subset \p^3$,
$p\notin H$, such that $f_1^{-1}(H\cap K')$,
if and only if $D$ is linearly equivalent to the
line section of $C$.

The curve $K'\subset\p^3$ as well as
the point $p$ and the plane $H$ are
uniquely determined up to a projective linear transformation of $\p^3$.
\end{prop} 
It is convenient to introduce homogeneous
coordinates $z_0:z_1:z_2:z_3$
to $\p^3$ so that $H$ is a horizontal
plane $\{z_3=0\}$ and $p=(0:0:0:1)$.

\begin{defn}\label{def-diagram}
Let $C\subset\p^2$
be a nodal real algebraic curve of degree $d-2$,
$\nu:K\to C$ is its normalization,
and $D=D^+-D^-$ (with disjoint effective divisors
$D^\pm$) be a real divisor on $K$
linearly equivalent to the pull-back of
the hyperplane section of $C$. 
%A pair $(C,D)$ consisting
%of an algebraic curve $\C C\subset\cp^2$
%of degree $d-2$ defined over $\R$, a divisor
%$\C Z\subset\C C$ of degree $d$
%and two distinct points $\{x'',y''\}\in\C C$.
%The triple $(\C C;\C Z,\{x'',y''\})$
We say that the pair $(C;D)$
is a {\em nodal diagram} if $\deg D^+=d$,
the divisor $D^-$ consists of two points of $K$
with distinct images under $\nu$,
%(so that $\deg D^+=d$),$
and the curve $K'\subset\p^3$ provided by Proposition
\ref{reconstruct-rkprime} has no singular points
other than $p$. 
\ignore{
\begin{itemize}
\item The curve $\C C\subset\cp^2$ is nodal,
i.e. all of its singularities are simple nodes.
\item The divisor $\C Z$ as well as the set
$\{x'',y''\}$ are invariant with respect to the
involution $\conj$ of complex conjugation.
\item The sum of the hyperplane divisor
(of degree $d-2$) of $\C C\subset\cp^2$ and
$\{x'',y''\}$ is linearly equivalent (over $\R$)
to $\C Z$.
\item The curve $\R K'\subset\rp^3$ obtained
from $(\C C;\C Z,\{x'',y''\})$ with the help
of Proposition \ref{reconstruct-rkprime}
does not have singularities other than the
double point $p\subset\R K'$.
\end{itemize}
}
\end{defn}
Note that since the only singularities
of $C=\tilde\pi_p(K)$ are its nodes,
the only possible singularities of 
$K'\setminus\{p\}$ are nodes of $C$ lifted by
$\tilde\pi_p$. But if $K'\setminus\{p\}$
is nonsingular then for each node $s\in C$
the two points of $\tilde\pi_p^{-1}$ are
distinguished by the value of the coordinate $\frac{z_3}{z_0}$.        

According to Proposition \ref{reconstruct-rkprime}
The pair $(C,D)$ determines the curve $K'\subset\p^3$.
Given  $D^-$ there might be several choices
of $D^+$ producing topologically equivalent
curves $K'$.
We may consider coarser data consisting only of 
$(C,D^-)$ and {\em topological data} of lifting
of $C$ to $\p^3$ as follows.

%Denote with $\Sigma=\R C\setminus\nu(\R K)$
%the (finite) set of nodes of $C$.
Consider the pair $(\R C,\R D^-)$ consisting
of a nodal algebraic curve
$\R C\subset\rp^2$
and the divisor
%$\{x'',y''\}
$\R D^-=D^-\cap\R K\subset K$ consisting of two distinct points
invariant under $\conj$ and such
that $\nu(\R D^-)$ is disjoint from
the set of nodes $\Sigma$.
%of elliptic nodes of $C$.
Let $\tau:\R K\to\rp^3$ be a smooth immersion
%lifting $\R C$ to $\rp^3$, i.e.
such that
$\tilde\pi_p\circ\tau=\nu$ and
$\tau|_{\R K\setminus\R D^-}$
is a proper embedding
of $\R K\setminus\R D^-$ to
$\rp^3\setminus\{p\}$.
Let $\sigma:\Sigma^e\to\{\pm 1\}$
be any function, where
$\Sigma^e\subset\Sigma$ is the set
of elliptic nodes. 

\begin{defn}\label{def-virtual}
An equivalence class of quadruples
$(\R C;\R D^-,\tau,\sigma)$
with respect to isotopies of $D^-$ and $\tau$
so that $D^-$ remains disjoint from 
the nodes of $C$ and
$\tau|_{\R K\setminus\tau^{-1}(p)}$
remains an embedding (the curve $C$ as well
as the function $\tau$ are fixed)
is called a {\em virtual nodal diagram} of $\R K'$.
%
%A {\em virtual diagram} is an equivalence
%class of triples $(\R C;D^-,\sigma)$
%with respect to isotopies of $D^-$ and $\sigma$
%within 
%enhanced
%with the identification between the branches
%of $\R\tilde K'$ and the inverse images under
%$\pi_p|_{S^2}$ near the special points of $\R C$ 
%is called {\em a virtual nodal diagram}.
\end{defn}
We depict virtual nodal diagrams in 
the same style as coventional knot diagram.
The points of $\R D^-$ are marked by bold
points. If $I\subset \R C$ is
a small open interval around a point from
$\R D^-$ then a half of this interval goes
very high up under $\tau$
(as the point itself goes to $p=(0:0:0:1)$)
while the other half goes very low down.
In our pictures we indicate the low half with a break (see e.g. Figure \ref{move-pole}).
\ignore{
\begin{rem}
The notion of virtual nodal diagram
%of a spatial algebraic curve
$\R K'$ is similar to that of diagram of a (smooth)
knot in $\R^3$: we record the information of
overcrossings and undercrossings. There are however
several differences. One is caused by our consideration
of $\rp^3$ as the ambient space, see e.g.
\cite{JuViro} for such set-up.
Another is the presence
of the points $x'',y''$ where the ``upper'' and
``lower'' branches of
$\R K'\setminus\{p\}$ are given by two ``half-arcs''.
Finally, our virtual nodal diagrams are defined
as algebraic curves.
Accordingly, we identify two virtual nodal diagrams
if their underlying real curves $\R C$
are {\em rigidly isotopic} in $\rp^2$
(i.e. isotopic in the class of nodal
real algebraic curves
of the same degree and genus).
%i.e.
%isotopic in the class of pairs $(\R C;\{x",y"\})$
%where $\R C$ is an algebraic curve of degree $d$
%in the plane
%with simple nodes disjoint from the points $x",y"$. 
\end{rem}
}

%To reconstruct $\R K'\subset\rp^3$
%from $(\R C;x'',y'')$ it suffices
%to consider a real divisor $Z$ of degree $d$ in $\C C$
%that is linearly equivalent to the union of a
%line section of the curve $\C C\subset\cp^2$
%and $\{x'',y''\}$.
%By a real divisor we mean a divisor invariant with
%respect to the complex conjugation in $\C C$.
%%The choice of $Z$ is not unique.
%The divisor $Z$ corresponds to a choice of 
%a plane $\R H\subset\rp^3$ disjoint from $p$.
%Clearly, its choice is not unique.

\begin{defn}\label{def-sep}
%Suppose that $\R C\subset\rp^2$
%is the (real algebraic) closure of the image $\pi_p(\R K'\setminus\{p\})$
%under the projection $\pi_p:\rp^3\setminus\{p\}\to\rp^2$.
We say that a plane $\R H\subset\rp^3$ is
%{\em determines the lift of $\R C$ to $\rp^3$}
{\em separating for $K'$}
%if $\R K'\setminus\{p\}\subset\rp^3$ is a smoothly embedded
%submanifold and
if for every line $l\subset\rp^3$
passing through $p$ and two distinct points
$u,v\in\R K'\setminus\{p\}$
the points $p$ and $\R H\cap l$
belong to different components of $l\setminus\{u,v\}$.
%$u\neq v$,
%the point $\R H\cap l_{u,v}$ is contained (strictly)
%inside the interval connecting $u$ and $v$ in
%$l_{u,v}\setminus\{p\}$.
\end{defn}
This means that the vertical coordinate function
$z_3/z_0$ have different signs
at $u$ and $v$.
%In other words $u$ and $v$ belong to different open
%half-spaces locally divided by $\R H$ when viewed
%from $p$.
Since any three points are coplanar we get
the following straightforward statement.
\begin{lem}\label{3pt-separate}
If the number of double points of $\R C$ is at most 3
then we may choose $\R H\subset\rp^3$ so that 
it is separating for $\R K'\subset\rp^3$.
\end{lem}

While in general it might be difficult to reconstruct
a virtual nodal diagram from a nodal diagram
there is a special case when it is easy. 
\begin{defn}\label{def3Dlift}
We say that
%a real divisor $\C D\subset\C C$
%(from the class of hyperplane section of $\C C$)
%{\em determines a 3D-lift}
the divisor $D$ is {\em 3D-explicit}
if the hyperplane
$\R H$ corresponding to $D^+$
is separating.
% (in the sense
%of Definition \ref{def-sep})
%for the curve $\R K'\subset\rp^3$
%defined by Proposition \ref{cd-rkprime}.
In such case we also call the nodal diagram 
$(C,D)$ 3D-explicit.
%$\Delta=(\C C;\C Z,\{x'',y''\})$ determines
%a 3D-lift if $\C Z\setminus\{x'',y''\}$
%determines a 3D-lift.
%This means that for every nodal point $s\in\R C$
%there exists a straight interval $I\subset\rp^3$
%such that $\dd I=\pi_p^{-1}(s)\cap\R K$,
%$I\cap\R H\neq\emptyset$ while $p\notin I$.
\end{defn}
%If $\C D$ determines a 3D-lift..

Since being 3D-explicit is determined by
the sign of the vertical coordinate, there
is the following straightforward criterion.
Let $Z\subset\R C$ be a 1-cycle, i.e.
the subspace of $\rp^2$ homeomorphic to a circle.
At a nodal point of $\R C$ the cycle $Z$
may stay on the same branch of $\R K$ or
may change it. We refer to the latter case
as the {\em corner} of $Z$. Let $n(Z)$ be
the number of the corners of $Z$ in the case
when $Z$ is null-homologous, and one
plus the number of corners in the case when
$Z$ is homologically non-trivial in $\rp^2$. 
\begin{prop}\label{prop-3Dliftcr}
%A real divisor $\C D=\C D_+-\C D_-$ with $\deg\C D_-=2$ on a real generically
%immersed curve $\C C$ determines a 3D-lift of $\R C$ 
Suppose that $\R C$ is connected.
The divisor $D$ is 3D-explicit if and only
if for each 1-cycle
$Z\subset\R C$
%(considered as a topological
%subspace of $\rp^2$ homeomorphic to a circle)
%null-homologous in $\rp^2$
the number of points from
$Z\cap D$ (counted with multiplicity) is congruent modulo 2
to $n(Z)$.
\end{prop}
\ignore{
\begin{proof}
We lift $Z\setminus\Sigma$ to $\rp^3$ so that 
the result $\tilde Z$ is contained in $\R K'$.
Then we connect the endpoints of $\tilde Z$ with the intervals in $\rp^3$
disjoint from the projection point $p$ to obtain a closed curve homologous
to zero in $\rp^3$. Such curve has even intersection number with any plane in $\rp^3$,
in particular with the horizontal plane $\R H$.
\end{proof}
}

Consider a deformation $\C D_t\subset\C C$,
$0\le t\le 1$,
of the divisor $\C D=\C D_0$ within
real divisors in the same
linear equivalence class (of the line section
of $\C C\subset\cp^2$). Assume that the degree
of the positive part $\C D^+_t$ 
(and thus also of the negative part $D^-$)
remains constant.

%%Recall that $\Sigma\subset\R C$ is the set of
%%the real nodes of $\R C$.
%Denote with $\tilde\Sigma=\nu^{-1}(\Sigma)\subset K$
%the inverse image of the nodal set $\Sigma\subset\R C$
%under the normalization.
%Denote with $\tilde q\subset\tilde\Sigma$
%the inverse image of $q\in\Sigma$. 
%%By Proposition \ref{reconstruct-rkprime}
%%$\C Z_t$ define a deformation $\R K'_t\subset\rp^3$
%%of the curve $\R K'$.
%%such that a
%%family of planes $\R H_t\subset\R P^3$
%%determine the lift of $\R C$ to $\rp^3$ for all $t$.
%%Clearly 
\begin{prop}\label{CD-ri}
Let $(C,D)$ be a 3D-explicit nodal diagram,
%$\C D=\C D_0$ determines a 3D-lift.
% lift of $\R C$
and $D_t$, $t\in [0,1]$, is a deformation
of $D=D_0$ in the class of real divisors
such that $(C,D_t)$ is a 3D-explicit
nodal diagram for $t\neq \frac12$
while for every $q\in\Sigma$ 
the following condition holds:
%Either
%both parts $\C D^+_t,\C D^-_t\subset\C C$
if $D_{\frac 12}^\pm\cap\nu^{-1}(q)=\emptyset$
then $D_{\frac 12}^\mp\cap\nu^{-1}(q)=\emptyset$.
%or one of the points of $\nu^{-1}(q)$
%one of the two points of $\nu^{-1}(q)$ is contained
%with weight $\pm 1$ in $\C D^+_t$ and the other with
%the opposite weight $\mp 1$ in $\C D^-_t$. 

Then
%for all $t\in [0,1]$ the pair $(C,D_t)$ is
%a nodal diagram,
%if $\C D=\C D_0$ determines a 3D-lift and
%$\C D_1$ is disjoint from $\tilde\Sigma$ then
%$\C D_1$ also determines a 3D-lift.
%while
the curves $K'_t$
%are smooth, so that they
provide a rigid isotopy between
$K'_1$ and $K'=K'_0$.
%Furthermore, $(C,D_t)$ is 3D-explicit for
%values of $t$ such that
%$\nu(D_t)$ is disjoint from $\Sigma$.
%defined
%by $\C D_1$ is rigidly isotopic to $\R K'=\R K'_1\subset\rp^3$.
\end{prop}
Here by a rigid isotopy we mean a deformation in the class of curves with 
a node at $\{p\}$
such that $K'_t\setminus\{p\}$ remains embedded.
\begin{proof}
%For the proof it suffices to note that the
%plane $\R H_t$ corresponding
%to $D^+_t$ 
If the support of $\C D_t$ is disjoint from
$\nu^{-1}(\Sigma)$ then the plane $\R H$
remains separating.
Otherwise the line $l_q$ in $\rp^3$ connecting
$p$ and a point in the image of
$\nu^{-1}(q)$ on $K'_t$, $q\in\Sigma$,
becomes tangent to $\R K'_t$
at $p$ while the intersection
$(\R K'_t\setminus\{p\})\cap l_q$
%in a unique point from
is contained in $\R H$. 
In both cases $p$ is the only singular point of $K'_t$, and it is a node.
\end{proof}

\newcommand{\wri}{\operatorname{wr}}
\newcommand{\self}{\operatorname{self}}
\newcommand{\wris}{\operatorname{wr}_{\operatorname{self}}}

\ignore{
%Suppose that the divisor $\R Z$ defines the lift to $\rp^3$.
We define the writhe $\wri(\R C,\C Z)$ as the sum
of the signs $\pm 1$ as in Figure \ref{writhe}
at the subset $\Sigma_{\self}\subset\Sigma\subset\R C$ that correspond 
to self-crossings of the same components of $\R K'$.
Since both branches of $\R C$ at
a point $q\in\Sigma_{\self}$ belong to the same
connected component, the sign is invariant of
the choice of orientation.
\begin{figure}[h]
\includegraphics[height=12mm]{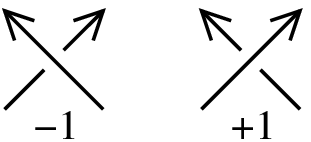}
\caption{Signs of diagram self-crossing points\label{writhe}}
\end{figure}
}

\section{Viro invariant through diagrams}
%We conclude this section with the computation of 
If $K\subset\p^3$ is a 
%$K_s$, $s=\pm 1$, be the result of positive 
%(resp. negative)
perturbation
of the nodal curve $K'\subset\p^3$
as in Lemma \ref{rkprime} then 
the Viro invariant $w(K)$
%and $\wl(K)$
can be computed
in terms of the virtual diagram $(\R C;\R D,\tau,\sigma)$.
%for the degeneration $\R K'\ni p$.
%Recall that we consider a real algebraic link
%$\R K\subset\rp^3$ obtained as a deformation
%(through a deformation of the projection of $\R \tilde K$
%obtained as a lift of $\R K$ to $\rp^4$)
%of an immersed curve $\R K'\subset\rp^3$
%with a single self-crossing point $p$.
Define $c=0$ if the two branches at $p$ 
belong to different components of $\R\tilde K$.
If both branches come from the same component of $\R\tilde K$ then
we orient $\R\tilde K$ arbitrarily and 
define $c$ as the sign of the double point
resulting from $p$
of the knot diagram of $\R K$ (when projected from
a point far from $p$).

Let $u\in\R C$ be a smooth point.
Choose a local orientation of $\rp^2$ near $u$
and an orientation of a component $M\subset\R C$
containing $u$. Note that it amounts to
a choice of generator in $H_1(\rp^2\setminus\{u\})$.

Let $u_+,u_-\notin\R C$
be points obtained by small deformations of $u$ to the
positive and negative side of $M$ respectively.
We define $\ind_M(u_\pm)$ as the image
of $M$ in $H_1(\rp^2\setminus \{u_\pm\})=\Z$
and set $\ind_M(u)=\frac{\ind_M(u_+)
+\ind_M(u_-)}2$. Clearly, the sign of this number
changes if we change the local orientation of $\rp^2$
or the orientation of $M$.
However, in the case when $u\in D^-$
%define the orientation of $\rp^3$
the orientation of $M$ defines the orientation of 
the tangent line to $p$ so that together with the  
local orientation of $\rp^2$ we can compare the
resulting orientation with the (standard)
orientation of
the ambient $\rp^3\supset\R K'$.
We set $i_M(u)\in\frac 12\Z$ to be $\ind_M(u)$
in these orientations agree and $-\ind_M(u)$
otherwise.

In the case when $K$
%(and thus $C$)
is a curve of type I
(see \cite{Rokhlin})
we can similarly define the index
$\ind_{\R C}(u)\in\frac 12
H_1(\rp^2\setminus \{u\})$
of $u$ with respect to the entire curve $\R C$ as well as the corresponding
half-integer number $i_{\R C}(u)\in\frac 12\Z$
using any of the two complex orientations of $\R C$.
Also in this case we define the linking number
$$\lambda(\R K)=\sum\limits_{M,N}\lk(M,N),$$
where the sum is taken over all pairs of different connected
components $M,N\subset\R K$ and the number $\lk(M,N)\in\frac 12\Z$
is the linking number in $\rp^3$ of the components $M$ and $N$
enhanced with orientations induced from a complex orientation of $\R K$.
As it was noted in \cite{Vi} for type I curves $K$
it is also useful to consider the invariant
$$\wl(K)=w(K)+\lambda(\R K).$$
In this case we define
%$c_\lambda$
%to be $c$
%if the two branches of $K'$ at $p$ 
%correspond to the same component of $\R K$.
%Otherwise we define
$c_\lambda=\pm1$ according
to the sign of resolution of $p\in\R K'$
with respect to the complex orientations of $\R K$,
so that we have $c_\lambda=\pm 1$ even if
the two branches of $\R K'$ at $p$ correspond
to different components of $\R K$.
Similarly, for a hyperbolic node $q\in\Sigma$
we define $\sigma_\lambda(q)$ to be the sign
of the corresponding crossing point with respect
to the complex orientation of $\R K$. 

\ignore{
If $K$ is of type I and $q\in\Sigma$
then in addition to the sign $\sigma(q)$
already defined we define the sign
$\tau(q)=\pm 1$ to coincide with $\sigma(q)$
(cf. Figure \ref{writhe}) unless $q$ is a hyperbolic
node corresponding to intersection of different components of the normalization
of $\R C$. In the latter case we define $\tau(q)$ 
also according to Figure \ref{writhe}
using any of the two complex orientations of $\R C$
(recall that in this case $\sigma(q)=0$ according to our convention).
}
\begin{prop}\label{viro-compute}
We have $$w=\sum\limits_{q\in\Sigma}\sigma(q)+
2\sum\limits_{u\in D^-}i_{M}(u)
+c.$$
Similarly,
$$\wl(\R K)=\sum\limits_{q\in\Sigma}\sigma_\lambda(q)+
2\sum\limits_{u\in D^-}i_{\R C}(u)
+c_\lambda$$
if $\R K$ is of type I.
\end{prop}
\begin{proof}
After scaling $\R^3=\rp^3\setminus\R H$ by
a very large number we may assume that $\R K$ is
obtained by a deformation of the union of $\R C$
with the two lines connecting the two points
of $D^-$ and $p$.
%The double points corresponding to the self-crossing
%of $\R C$ as well as to the intersection of 
%the two lines are deformed by moving the branches apart,
%so they contribute $\wris(\R C,Z)+c$ to $w$.
The points of intersections of these lines with $\R C$
get smoothed.
The remaining intersection points contribute
$2\sum\limits_{u\in D^-}i_{M}(u)$
to $w$
and
$2\sum\limits_{u\in D^-}i_{\R C}(u)$
to $\wl(\R K)$.
\end{proof}

%% file: upto5.tex
\subsection{Links of degree 4 and lower}
Let us apply Lemma \ref{lRR} and Corollary \ref{dle6}
for smooth irreducible algebraic curves $\R K\subset\rp^3$ of small degree $d$. 
%All smooth irreducible degree one and degree two algebraic curves in $\rp^3$
If $d=1,2$ then $g=0$ and $\R K$ is (algebraically) planar. If $d=3$ and $g=1$ then
$\R K$ is also planar. If $d=3$ and $g=0$ then  
$\R K$ is hyperboloidal of bidegree $(2,1)$, and thus
it is only topologically planar.

Consider the case of $d=4$. 
By Corollary \ref{dle6} any such link sits on a quadric.
If $g=3$ then it is a planar link.
%$p_0^l$ with $0\le l\le 4$.
Otherwise $\R K$ corresponds
to a bidegree $(a,b)$ curve with $a+b=4$ and $(a-1)(b-1)=g$.
Therefore we never encounter $d=4$, $g=2$ curves.
If $g=1$ then $(a,b)=(2,2)$ and thus $L$ must be
topologically planar
by Proposition \ref{plan-hyp}.
Finally, if $g=0$ we have $(a,b)=(3,1)$ and thus $\R K$ 
is of hyperboloidal type $h_{3,1}$.

\subsection{Links of degree 5}
By Corollary \ref{dle6} if $d=5$ and $g>1$ then $\R K$ is hyperboloidal.
If $g=6$ it is planar. The bidegree of a smooth irreducible curve 
$\R K\subset\R Q\subset\rp^3$ is either $(3,2)$ or $(4,1)$.
In the first case we have $g=2$ with a topologically
planar
%type $p_1^l$, $l=1,2,3$ of the
link by Proposition \ref{plan-hyp}.
In the second case we have $g=0$.
Therefore, the cases $d=5$ and $g=3,4$ never appear.

\begin{figure}[h]
\includegraphics[width=100mm]{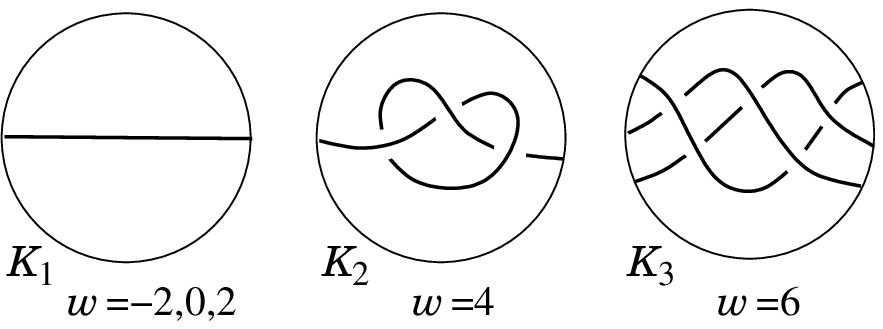}
\caption{Rational quintic knots \label{d5g0}}
\end{figure}

The following result was obtained by Bjorklund \cite{Bj}.
Recall that the topological isotopy type of $\R K$ is the
equivalence class of
the pair $(\rp^3,\R K)$ up to homeomorphism.
Note that an orientation-reversing homeomorphism
takes $w(K)$ to $-w(K)$.
Recall that we say that that %$\R K_1$ and $\R K_2$
two real algebraic curves embedded in $\p^3$
are rigidly isotopic if one can be deformed
to another in the class of embedded smooth curves.
\begin{thm}[Bjorklund \cite{Bj}]
\label{thm-bj}
There are three distinct topological isotopy types of
$(\rp^3,\R K)$ for $d=5$, $g=0$
shown at Figure \ref{d5g0}.
\begin{itemize}
\item
The trivial knot $K_1$.
% is of the planar isotopy type $p_1$ (that of $(\rp^3,\rp^1)$). 
In this case $w=\pm 2$ or $w=0$.
\item
The long trefoil knot $K_2$
(a connected sum of a trefoil and a projective line).
%, see Figure \ref{longtrefoil}).
In this case
$w=\pm 4$.
\item
The hyperboloidal knot $K_3$
of type $h_{4,1}=\bar T(6,4)$.
In this case $w=\pm 6$.
\end{itemize}
Furthermore, any two smooth curves of degree 5
and genus 0 in $\rp^3$ are rigidly isotopic if
and only if they have the same invariant $w$.
\end{thm}

Let us pass to the case when $\R K\subset\rp^3$ is
a (non-empty) smooth degree $d=5$, genus $g=1$ curve.
As the number of components of $\R K$ is not greater
than $g+1$ by Harnack's inequality \cite{Harnack}
we have two possibilities:
either $\R K$ is connected,
or it contains two connected components.
In the latter case the quotient $\C K/\conj$ of $\C K$ by the involution
$\conj$ of complex conjugation is an annulus and thus $\R K$
must be of type I. 
Thus the invariant $\wl$
is well-defined.
In the former case the quotient $\C K/\conj$ is a M\"obius band,
so that $\CC K\setminus\R K$ is connected, i.e. $\R K$ is of type II, so we are restricted to consideration of $w$.

\begin{figure}[h]
\includegraphics[width=100mm]{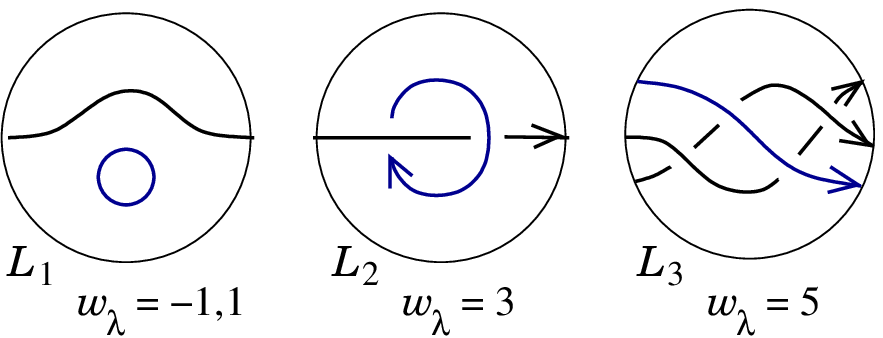}
\caption{Elliptic quintic two-component links. \label{d5g1}}
\end{figure}

\begin{thm}
\label{thm-d5g1}
There are three distinct topological isotopy types of $(\rp^3,\R K)$
for $d=5$, $g=1$, in the case when $\R K$ is a two-component link, see Figure \ref{d5g1}.
\begin{itemize}
\item
%The curve of type I %$\R K\subset\rp^3$
%consisting of two components:
%a line $\rp^1\subset\rp^3$
%and an unknotted circle around this line
%(so that the circle bounds an embedded disk
%intersecting the line transversely in a single point).
%In this case we have
The trivial (planar) link $L_1$.
In this case
$w=\pm 1$, $\wl=\pm 1$.

\item
The link $L_2$
%The curve of type I %$\R K\subset\rp^3$
consisting of a line $\rp^1\subset\rp^3$
and an unknotted circle around this line.
%(so that the circle bounds an embedded disk
%intersecting the line transversely in a single point).
%In this case we have
In this case
$w=\pm 1$, $\wl=\pm 3$.\\
(Figure \ref{d5g1} shows a complex orientation 
in the case $w=1$.)

\item
The link $L_3$ consisting of a hyperboloidal 
knot of type $h_{3,1}$ and a line $\rp^1\subset\rp^3$
disjoint from the hyperboloid containing the other components. 
In this case
$w=\pm 3$, $\wl=\pm 5$.
\end{itemize}

If $d=5$, $g=1$ and $\R K$ is connected 
then it is isotopic to 
$\rp^1\subset\rp^3$ (see $K_1$ from Figure \ref{d5g0}).
In this case we have
$w=\pm 1$.
%\end{itemize}

Furthermore, all two-component real algebraic links of degree 5
and genus 1 in $\rp^3$ with the same value of $\wl$ are rigidly isotopic.
%are rigidly isotopic if
%and only if they have the same value of $\wl$.
Also all connected real algebraic knots of degree 5
and genus 1 in $\rp^3$ with the same value of $w$
are rigidly isotopic.
\end{thm}
%Note, in particular, that 
%the value of $\wl$ (when it is defined)
%for a two-component link
%of degree 5 and genus 1 determines the value of $w$.
%If $v_\rho+2\lambda=5$ then $v_\rho=3$ and $\lambda=1$. 
%If $v_\rho+2\lambda=3$ then $v_\rho=1$ and $\lambda=1$. 
%If $v_\rho+2\lambda=5$ then $v_\rho=1$ and $\lambda=0$. 
\begin{proof}
The rank $r$ of the linear system defined by
the plane section of $K\subset\p^3$ is at least
$d-g=4$
by the Riemann-Roch theorem. By Lemma \ref{rkprime}
we may assume that $K$ is obtained by
deformation of a curve $\R K'$ 
with a self crossing point $p$. 
By Proposition \ref{pirkprime} %$\pi_p(\R K')$
%the curve $C$ in the corresponding nodal diagram
%$(C,D)$
%%$\R C\subset\rp^2$ of $\R K'$
%is a nodal cubic of a genus 1.
%% curve of degree 3 which is necessarily an embedded
%By the adjunction formula $C\subset\p^2$ is smooth.
%
%In other words
the curve 
%rigid isotopy type of $K\subset\p^3$
%is determined by a smooth planar cubic curve
$C=\tilde\pi_p(K')\subset\rp^2$
in the corresponding nodal diagram is cubic of genus 1.
Thus $C$ is smooth.
%enhanced with
%the choice of two distinct points $x'',y''\in \C C$
%invariant with respect to the involution $\conj$
%of complex conjugation.
%The equivalence class of a divisor $D$
%with $\deg D^-=2$ is determined by the distribution

Suppose that $D^-\cap\R C=\emptyset$.
Note that this determines the equivalence class
of $D$ up to real deformations.
% form a complex conjugate
%pair, i.e. $x'',y''\notin\R C$.
In this case $\R K'$ is topologically isotopic
to the union of the planar cubic curve
isotopic to $\R C$
%$\R C\subset\rp^2\subset\rp^3$
and a solitary node at $p$.
After a deformation the solitary node $p$ disappears.
%Here the auxiliary $\rp^2\subset\rp^3$
%is assumed to be disjoint from $p$.
By Proposition \ref{viro-compute} we
have $w(K)=c=\pm 1$.
%while both signs of $c$ appear by
%Lemma \ref{rkprime}.

%The image under $\pi_p$ 
%of a perturbation $\R K$ of $\R K'$ thus has an elliptic
%node and two pairs of complex conjugate nodes. 
%Thus we get the two last cases of the theorem with
%$v_\rho=\pm 1$. Both signs of $\pm 1$ appear by
%Lemma \ref{rkprime}.

In other cases we have $D^-\subset\R C$.
Then $\R K'$ is obtained from $\R C\subset\rp^2\subset
\rp^3$
by attaching the two lines connecting $p$ with
the points of $D^-$ and %$x''$ and $y''$
then perturbing the result with the help of $D^+$.
Note that up to equivalence $D^+$ is determined
by the parity of the number of points in each
connected component of $\R C\setminus D^-$. 

%can be deformed one to another and thus are 
%equivalent.
%By Proposition \ref{reconstruct-rkprime}
%the perturbation is determined
%by a divisor of degree 5 linearly equivalent to 
%the union of the plane section with $x''$ and $y''$.
%To get $\R K$ we resolve the self-crossing at $p$
%in each of the two possible ways.

Suppose that 
%the smooth cubic curve $\R C\subset\rp^2$ is connected
%both $x''$ and $y''$ belong to
$D^-$ is contained in the homologically
non-trivial component $J$ of $\R C\subset\rp^2$.
Since $D$ is linearly equivalent to the hyperplane
section of $\R C$ we must have even number of points in
$\R C\setminus J$ and different parities
in the two arcs of $J\setminus D^-$.
%it is of type I and then 
%the corresponding component of $\R K$ %is isotopic
%to $\rp^1\subset\rp^3$ and 
Thus all such choices of $D$ are equivalent.
Once again, $\R K\subset\rp^3$ is topologically isotopic
to a curve sitting in a plane and isotopic to $\R C$.
%does not depend on the type of perturbation
We have $w(K)=c=\pm 1$.
%This is the only possibility 

If $\R K$ is of type II
(i.e. $\R C$ is connected) then $J=\R C$
and there are no other possibilities for $D$.
If $\R C$ is of type I (i.e. $\R K$ is a two-component
link) then $\wl(K)$ in both cases considered above
coincides with $w(K)$.

%If $\R K$ is a two-component link then
%$\R C$ consists of a non-trivial component and an oval,
%so that one or both points $x'', y''$ may sit on an oval.
If $\R K$ is a two-component link we also have
additional cases.
If $D^-$ has points in different components
of $\R K$ this again determines the class of equivalence
of $D$. Then $\R K$ is topologically isotopic
either to $L_1$ or $L_2$
%sits on the oval then we are either in
%the second or the third case of the theorem 
depending on the resolution at $p$.
We have $\wl(K)=\pm 2 + c_\lambda$ in these cases
by Proposition \ref{viro-compute}.

If $D^-\subset \R C\setminus J$ 
%both points sit on the oval
then there are two equivalence classes of $D$.
In one case we have odd number of points of $D^+$
in all three components of $\R C\setminus D^-$.
Then $\wl(K)=\pm 4 + c_\lambda$
and the topological type is $L_2$ or $L_3$
accordingly.
%we are either 
%in the first or the second case of the theorem.
%The value of $w$ is determined by Proposition
%\ref{viro-compute}.
In the other case both arcs of $\R C\setminus (J\cup D^-)$
have even number of points from $D^+$.
Then  $\wl(K)=c_\lambda$
and the topological type is $L_1$.
\ignore{
As $\R C\subset\rp^2$ is embedded,
%Note that all these cases can be 
%realized with the help of
an arbitrary divisor
$\C Z$ of degree 5 disjoint from $x''$ and $y''$
yields a nodal diagram, so all the cases above
are realizable.  
This completes the topological classification of 
the spatial quintic curves of genus 1.
}

To deduce the rigid isotopy classification
it remains to prove that in our construction
above the curves with the same $\wl$
that are obtained from different cases
considered above.

If $\wl=3$ then there are two options for the distribution
of $D^-$ between the components of $\R C$.
If $D^-\subset \R C\setminus J$ and $c=-1$
then %correspond to a self-intersection
$p$ corresponds to self-crossing 
of the topologically trivial (even) component of $\R K$.
If $D^-\cap J\neq\emptyset$,
$D^-\cap (\R C\setminus J)=\emptyset$ and $c=+1$
then $p$ corresponds to the intersection point
between different components of $\R K$.

If $\wl=1$ then there are three options for the distribution
of $D^-$ between the components of $\R C$.
They correspond to self-intersection of an even
%(null-homologous in $\rp^3$)
component of $\R K$,
%(the null-homologous component in $\rp^3),
the self-intersection
of an odd component of $\R K$
and the intersection point of distinct components
of $\R K$.
%do appear all.

We claim that in the case $|\wl|\le 3$ there exists
a real deformation
of $K$ to an immersed curve $K'$ with
a single crossing point corresponding
to distinct components of $K$.
To see this we consider the projection of $K$
from a generic point $p$ on the even component of $\R K$.
The image $B\subset\p^2$ of the projection is a quartic curve
with two odd connected components of the normalization
$\R K$.
Being odd, these components must intersect at a point $r\in\rp^2$.
Since the curve $B$
is elliptic, there is a second nodal point of $\R B$ which must be
either a self-intersection point $s\in\rp^2$ of a component $P\subset\R K$ 
or an elliptic double point $s\in\rp^2$ (the intersection of
two complex conjugate branches of $\C B$).

If $s$ is elliptic then it corresponds to
a pair $P_s\subset\CC K\setminus\R K$
of complex conjugate points in $\R K$
while $r$ corresponds to a pair $P_r\subset\R K$
of points
from different components of $\R B$.
Choose a plane $\R H$ passing through $P_s$
%a pair  of points corresponding to $s$
and a point $p^1_r\subset\ P_r$
from the odd component of $\R K$.
%The resulting divisor is equivalent to $D
The divisor $\C H\cap\C K$ on $K$ is
equivalent to $D^+$ and
consists of $P_s$, $p^1_r$
as well as another pair $P_m$ of points on $K$. 
If $P_m$
is contained in the odd component  of $\R K$
then we can deform $P_m$ into $\CC K\setminus\R K$
not changing the linear equivalence class.
If $P_m\cap\R K=\emptyset$ then we can further
deform $P_m$ to the even component of $\R K$
within the same linear equivalence class.

If $P_m$ is contained in the even component of $\R K$
then we can deform $P_m$ (moving the image
$p_B\in\R B$ of the projection point $p\in\rp^3$
if needed) so that $D^+-\{p_B\}$
stays in the same linear equivalence class
while the result of deformation of $D^+$
contains $P_r$. In this case the spatial curve
corresponding to $(C,D)$ by Proposition
%\ref{cd-rkprime}
\ref{reconstruct-rkprime}
has a double point at $r$.
 
If $s$ is not elliptic then it corresponds to
a self-intersection of one of the components
of $\R K$. We choose a plane $\R H\subset\rp^3$
passing through a point $p_r\in\R K$
and so that it separates
the pair $P_s$ in the sense of
Definition \ref{def-sep}.
%(note that now $P_s\subset\R K$).
Here we choose $p_r$ to be on the component
$A\subset\R K$ containing the pair $P_s$.

Recall that we assume that $|\wl|\le 3$.
By Proposition \ref{viro-compute}
this implies that if $A\cap D^+$ consists of
more than 4 points then two of them bound
an open interval $I\subset A$ disjoint from
$D^+$ and $p_B$. Thus a pair of points of $D^+$
can be pushed to $\CC K\setminus\R K$ and then 
to the other component $A'\subset\R K$.
Note that $\tilde\pi_p|_{A'}$ is an embedding
since $P_s\subset A$.
Thus we ensure that $D^+\cap A'$
consists at least of two points.
As in the case when $s$ is elliptic we deform
these points (along with $p_B$ if needed) to
ensure a crossing point between two different
components of the spatial curve.

Thus any embedded curve with $d=5$, $g=1$
and $|\wl|\le 3$ is obtained
by perturbing a nodal spatial curve with
a node corresponding to crossing of different
components. Therefore $\wl$ determines the rigid
isotopy type of type I $d=5$, $g=1$ real algebraic
links.
\end{proof}

%% file: deg6.tex
\section{Rational knots of degree 6}
\begin{figure}[h]
\includegraphics[width=100mm]{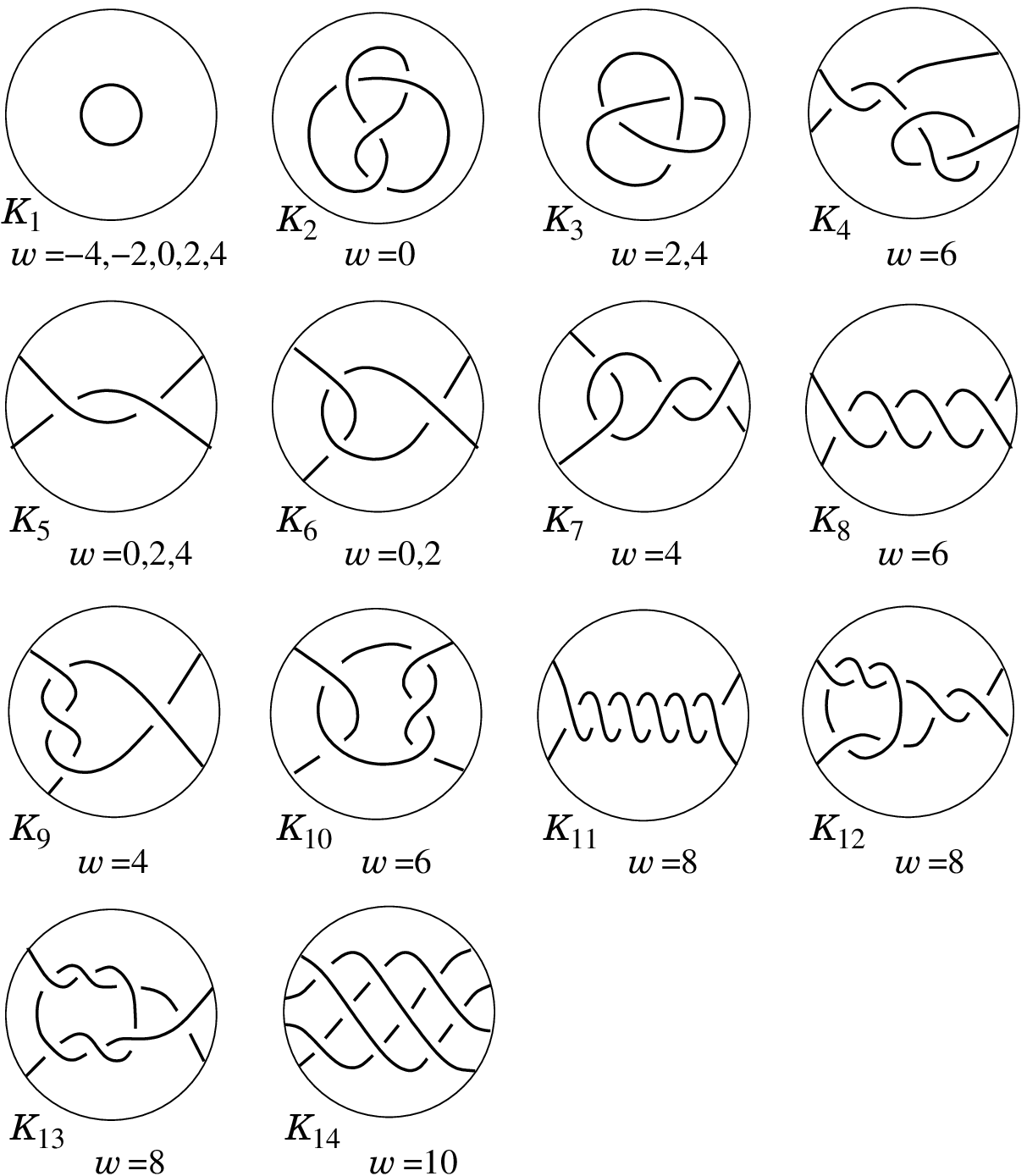}
\caption{
Real algebraic knots of degree 6 and genus 0.
\label{deg6c}}
\end{figure}

%\subsection{Statement of the theorems
%of topological and rigid classifications}
%For $d=6$ we start our classification from the case $g=0$.
\begin{thm}\label{thm-g0t}
There are 14 topological isotopy types (homeomorphism classes
of pairs $(\rp^3,\R K)$)
of rational real algebraic curves $\R K$ of
degree 6 embedded in the projective space $\rp^3$.
% considered up to
%diffeomorphism of the pair $(\rp^3,\R K)$.
Figure \ref{deg6c} lists the knots.
\end{thm}

\newcommand{\ew}{\operatorname{ew}}
\begin{thm}\label{thm-g0r}
There are 38 rigid isotopy types of rational real algebraic curves of degree 6
embedded in the projective space $\rp^3$.
Namely, each curve depicted on Figure \ref{deg6c} enhanced with 
a choice of the listed value for $w$ gives rise to
one rigid isotopy class of $\R K$ in the depicted knot type. 
Furthermore, simultaneous reflection of $\R K$ in $\rp^3$ and changing
the sign of $w$ gives a new rigid isotopy type with the exception
when the knot type of $\R K$ is amphichiral (the first two 
knots of Figure \ref{deg6c}).
\end{thm}

\subsection{Quartic nodal diagrams and odd arcs}
Lemma \ref{rkprime} and
Proposition \ref{reconstruct-rkprime}
%and \ref{cd-rkprime}
%to a smooth algebraic curve $\R K\subset\rp^3$ of degree 6 we
reduce Theorems \ref{thm-g0t} and \ref{thm-g0r}
to classification
of the nodal diagrams %$(\C C;\C Z,\{x'',y''\})$
$(C,D)$ with respect to equivalences corresponding to 
topological and rigid isotopies of the resulting
spatial curves.

Since $d=6$ the curve $\C C$
is a nodal quartic. Thus it has at most three
real nodes. 
%of generically immersed quartics $\R C=\pi_p(\R K')\subset\rp^2$ enhanced with
%a real divisor $\C D=\C D_+-\C D_-$ of degree 4 and $\C D_-=\{x'',y''\}$
%consisting of two points distinct from the self-crossings of $\C C$.
%Since $\R C\subset\rp^2$ is a generically immersed quartic curve, the 
%cardinality of the set $\Sigma$ 
%of its real nodes is at most 3.
By Lemma \ref{3pt-separate} we may assume that
$\C D$ is 3D-explicit.
%Thus the diagram
%$(C,D)$
%$(\C C;\C Z,\{x'',y''\})$
%is subject to
Thus we may apply Proposition \ref{prop-3Dliftcr}.

Given a nodal diagram of $\R K'$
we may consider
parity of the number of elements of $A\cap\C D_+$
(counted with multiplicities)
for any connected component
$A\subset\R C\setminus(\Sigma\cup\C D_-)$.
%Any such component lifts to a proper arc connecting $\R^3\setminus\R H$, so its intersection
%with $\R H$ is even if and only if the endpoints are from the same half-space
%of $\R^3\setminus\R H$.
%Namely, for an arc intersecting the infinite line $\rp^1\subset\rp^2$...
We refer to such components $A$ as {\em  diagram arc}.
A diagram arc $A$ is called {\em odd} if the parity of $A\cap\C D_+$
is odd, and {\em even} otherwise.
Clearly, the parity of a diagram arc depends only
on the underlying virtual knot diagram,
so that we may speak of odd arcs of virtual
nodal diagrams.
%Namely, a diagram arc is odd if the identification
%between

\ignore{
\begin{defn}
A diagram consisting of a generically immersed 
real curve $\R C\subset \rp^2$ with the nodal set $\Sigma$,
two distinct points $x'',y''\in\R C\setminus\Sigma$ as
well as
the upper/lower half-space information near $x'',y''$ and the nodal set $\Sigma\subset \R C$
is called an {\em enhanced diagram} of $\R C$.
\end{defn}
Clearly, the notion of an odd diagram component makes sense for any
enhanced diagram of $\R C$ (and not only for $\R C=\pi_p(\R K')$).
}

\begin{prop}\label{prop-diagram}
Suppose that $(\C C;D^-,\tau,\sigma)$
is a virtual nodal diagram such that
all nodes of a rational quartic curve $\R C$
are hyperbolic.
The diagram $(\C C;D^-,\tau,\sigma)$
is realizable by a $3D$-explicit
nodal diagram $(C,D)$
%of
%a rational real curve $\R K'\subset\rp^3$ of degree 6
%with a double point at $p\in\rp^3$
if and only if the number of its odd arcs
is at most 6. 
%In this case a nodal diagram $(C,D)$
%for
%such a realizable $(\C C;\{x'',y''\})$
%may be chosen to be 3D-explicit.

Furthermore, in this case the virtual diagram $(\C C;D^-,\tau,\sigma)$ and the sign $c$
of deformation of $K'$
determine the embedded real algebraic curve $K\subset\p^3$ up to rigid isotopy.
%any two real rational curves of degree 6 
%obtained from curves corresponding
%to the same virtual nodal
%diagram with with nodal diagrams determining a 3D-lift
%corresponding
%to $(\C C;D^-,\tau,\sigma)$ with the same value of $c$
%the same virtual nodal diagram are 
%are rigidly isotopic.
%isotopic in the class of rational real curves with
%nodal diagrams determining a 3D-lift corresponding
%to the same virtual nodal diagram.
\end{prop} 
\begin{proof}
If the number of odd diagram arcs
is greater than 6 then so must be the degree
of the effective divisor $\C D^+$.
Conversely, if this number is at most 6 then we
may construct $\C D^+$ by selecting a point at each odd arc.
If needed we add pairs of conjugate points on $\C C$ to ensure $\deg\C D^+=6$.
The space of these choices is connected as any
pair of points of $D^+$ on the same arc may be deformed
off $\R C$ into $\CC C\setminus\R C$.
\end{proof}

\subsection{Moves of the diagrams}
%We restrict our attention to
%nodal diagrams determining a 3D-lift.
%Proposition \ref{prop-diagram} reduces Theorem \ref{thm-g0t} to
%classification of the knots coming from resolution of the double
%point of the spatial curves $\R K'$ corresponding to
%the virtual nodal diagrams of quartics with not
%more than 6 odd diagram components.
%Meanwhile some different
%enhanced diagrams correspond to the same spatial curve $\R K'$. 

%\begin{defn}
%A real divisor $\C D$ on $\C C$ is called {\em $\R$-mobile} if 
%\end{defn} 
Some nodal diagrams $(C,D)$ with distinct $C$ and $D$ correspond
to the same rigid isotopy class of knots.
We formulate the following
straightforward proposition
for nodal diagrams of rational curves
of arbitrary %node diagrams and
degree $d$.
%or $g=0$.
%In particular we shall apply in
%for the case of $d=6$, $g=1$ later in the text.

\begin{prop}\label{prop-moves}
%Let $K_1,K_2\subset\rp^3$ be
%two smooth algebraic curves 
%obtained from the curves $K'_1,K'_2\subset\rp^3$
%corresponding to two 
Suppose that nodal diagrams $(C_1,D_1)$
and $(C_2,D_2)$ of degree $d$ 
on rational curves $C_j$, $j=1,2$, of degree $d-2$
%is obtained from
%another nodal diagram $(C_2,D_2)$
are 3D-explicit and
related by means of one of the moves listed below.
%while
%with the same
%via resolving the (only) nodes of $\R K'_1$
%and $\R K'_2$ according to the same sign.
%Suppose in addition that the curves
%$\R C_1$ and $\R C_2$ from $\Delta_1$ and $\Delta_2$
%are topologically isotopic if
%$\Delta'$ is obtained from $\Delta$ by means of one of the following moves.
%(Moving a pole as well as the first Reidemeister moves
%are depicted on Figure \ref{moves}).
%Suppose that $\Delta_2$ is obtained from
%$\Delta_1$
%(with the underlying curve $\R C_1\subset\rp^2$
%enhanced with a pair of smooth points
%$x_1'',y_1''\in\R C$)
%by means of one of the following 5 moves
%and one of the diagrams $\Delta_j$, $j=1,2$,
%determines a 3D-lift.
%Assume that %one of the diagrams
%$(C_1,D_1)$ %, $j=1,2$,
%is 3D-explicit.

%Then %both diagrams $(C_j,D_j)$, $j=1,2$, are
%$(C_2,D_2_$ is also 3D-explicit.
%Furthermore, %in the case of $g=0$
Then the embedded real algebraic curves
$K_1, K_2\subset\p^3$
obtained from the corresponding curves
$K'_1$ and $K'_2$ by resolving
their double points in coherent directions
are rigidly isotopic.
% if $g=0$.
%Furthermore, if $g=1$ and the diagrams
%$(C_j,D_j)$, $j=1,2$,
%have at most $d-2$ odd arcs
%the curves $K_1,K_2\subset\p^3$ are also rigidly isotopic.
%\begin{itemize}
%\item

\noindent $\bullet$ {\bf (Moving a pole)}  {\em
Suppose that $C_1=C_2$.
Let $u\in\Sigma$ be a hyperbolic node of
$\R C=\R C_1=\R C_2$
and $u^+,u^-\in\R K$ be two points corresponding
to $u$ under the normalization $\nu:K\to C$.
Let $u^{\pm}_1$ be the result of moving $u^\pm$
a little along $\R K$ in one (arbitrarily chosen)
direction, and $u^{\pm}_2$ be the result of 
moving $u^{\pm}$ in the opposite direction. 
Let $D_+$ be an arbitrary real
effective divisor on $K$ of degree $d-1$ disjoint
from $\nu^{-1}(\Sigma)$, and $D^-$ be
an arbitrary real effective
divisor on $K$ disjoint from
$\nu^{-1}(\Sigma)$ and $D^+$. 
Define $D^{\pm}_j=D^{\pm}\cup\{u^{\pm}_j\}$, $j=1,2$. 
%$\R C_1\setminus D_1$
%contains an arc $A$
%connecting %$x''$ (or $y''$)
%a point $x_1\in D^-_1$ to a point $u\in\Sigma$.
%Suppose that $\R C_1\setminus D_1$ also contains 
%an arc $B$ connecting $u$ to a point $y_1\in D^+_1$ 
%Define $x_2\in\R C$ %D_2\setminus x_1$ 
%to be the point obtained
%diagram on the same curve
%$\R C_2=\R C_1$ by moving $x_1''$ (resp. $y_1''$)
%by moving $x_1$ along $A$ past $u$,
%and $y_2$ to be the point obtained by moving 
%$y_1$  along $B$ past $u$.
%and reversing the
%identification information between the
%branches of $\R\tilde K_1'$
%(see Definition \ref{def-virtual}) at $u$,
%while keeping the rest of the diagram $\Delta_1$ unchanged.
%Define $D^-_2=(D^-_1\setminus\{x_1\})\cup \{x_2\}$ and
%$D^+_2=(D^+_1\setminus\{y_1\})\cup \{y_2\}$.
(See Figure \ref{move-pole} for the change of
the corresponding virtual nodal diagram.)
}
\begin{figure}[h]
\includegraphics[height=12mm]{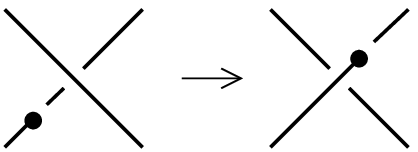}
\caption{Moving a pole. \label{move-pole}}
\end{figure}

\noindent $\bullet$ {\bf (Annihilation of two poles)}
{\em Suppose that $C_1=C_2$,
and $D^+$ is an arbitrary real effective divisor
of degree $d$ disjoint from $\Sigma$.
Let $x\in\R K\setminus(\nu^{-1}(\Sigma)\cup D^+)$
be a point.
Define $D^-_1$ to consist of two distinct points  
in $\R K$ close to $x$ and $D^-_2$ to be a pair
of conjugate points in $\CC K\setminus\R K$ close to $x$.
%D_1)$  contains an even diagram arc $A$ containing
%$x_1''$ and $y_1''$. Define $\Delta_2$ to be the diagram on the same
%curve $\R C_1=\R C_2$ by moving $x_1''$ and $y_1''$
%towards each other and away from the real diagram
%to a pair of conjugate points in $\C C_2\setminus\R C_2$
}

The remaining three moves are algebro-geometric
counterparts
of the Reidemeister moves from knot theory.
Here $C_1\neq C_2$, instead they are obtained
as different real perturbations
(within the class of real algebraic curves of degree $d$)
of a certain curve $C_0\subset\p^2$.
We denote the normalization of $C_0$ with
$\nu_0:K_0\to C_0$
and consider a certain divisor $D_0=D^+_0-D^-_0$
on $K$ obtained as the difference of
disjoint effective divisors $D^{\pm}_0$,
$\deg D^+_0=d$, $\deg D^-_0=2$, with $D^{\pm}_0$
disjoint from the normalization of the nodes of $C_0$. 
In the following moves the divisor $D_j$, $j=1,2$,
on the normalization $K_j\to C_j$ is obtained
as an arbitrary small deformation of the divisor
$D_0$ on $K_0$. 

\noindent $\bullet$ {\bf (Reidemeister 1)} {\em
%Suppose
%that $C_1$ and $C_2$ are obtained by
%two different real perturbations
%(within the class of real algebraic curves of degree $d$)
%of a curve 
%that
Here $C_0$ is a curve with a single cusp $u$
and simple nodes as all other singularities,
%Let $D_0=D^+_0-D^-_0$, $\deg D^+_0=d$,
%$\deg D^-_0=2$, be the difference of
%disjoint effective divisors on the normalization
while $D_0$ is such that $D^+_0\ni u$.
%We suppose that $\R C_1$ perturbs the cusp of $\R C_0$
%to a hyperbolic node %$\Lambda$
%disjoint from $x_1''$ and $y_1''$ in $\Delta_1$
%while $\R C_2$ perturbs the cusp of $\R C_0$
%to an elliptic node. %real point.
%Define $D_1$ to be the divisor on $K_1$
%obtained as a small perturbation of the divisor $D_0$
%such that $u\in D^+_0$ gets perturbed inside
%the loop $\Lambda$. Define $D_2$ on $K_2$
%to be an arbitrary small perturbation of $D_0$.
%$ to coincide with $\Delta_1$ everywhere away from
%the cusp.
See Figure \ref{reidem1} for the corresponding virtual
nodal diagrams.
}
%\begin{figure}[h]
%\includegraphics[height=12mm]{risunkiStepa/reidem1.eps}
%\caption{First Reidemeister move and a real solitary point. \label{reidem1}}
%\end{figure}

\noindent $\bullet$ {\bf (Reidemeister 2)} {\em
%Suppose that
%$C_1$ and $C_2$ are obtained by
%two different real perturbations
%(within the class of real algebraic curves of the same degree and genus)
%of a curve
Here $C_0$ has a single tacnode $u$
and simple nodes as all other singularities,
while the divisors $D^{\pm}_0$ are disjoint from $u$.  
%We suppose that $C_1$ perturbs the tacnode to two real
%nodes %connected by two small even arcs of $\R C_1$ disjoint from $x_1''$ and $y_1''$,
%while $C_2$ perturbs the tacnode into
%two elliptic nodes.
%solitary real points. Define  $\Delta_2$ on $\R C_2$ to coincide with $\Delta$
%everywhere away from the tacnode.
%Define $D_j$, $j=1,2$, to be arbitrary small
%deformations of $D_0$.
}

\noindent $\bullet$ {\bf (Reidemeister 3)} {\em
%Suppose that $\R C_1$ and $\R C_2$ are obtained by
%two different real perturbations
%(within the class of real algebraic curves of the same degree and genus)
%of 
Here $C_0$ is a curve with
a single ordinary triple point $u$
%with non-tangent branches
and simple nodes as all other singularities,
while $D^+_0$
contains a single point in the set $\nu_0^{-1}(u)$
(of cardinality 3).
%is compatible with the third Reidemeister move, i.e. the three
%small arcs connecting the nodes of $\R C_1$ resulting from the triple point
%are disjoint from $x_1''$ and $y_1''$ while at least of these arcs
%is even.
%Define $\Delta'$ according to the third Reidemeister move.
}
%\end{itemize}
\end{prop}
\begin{coro}\label{coro-2ptg0}
Any embedded real rational curve of degree 6 is rigidly isotopic
to a curve obtained from the curves
whose virtual nodal diagrams
are listed on Figure \ref{g0-2pt}
by resolving them according to $c=\pm 1$ 
as in Lemma \ref{rkprime}.

Each diagram $N^c_\epsilon$,
where $N$ is the number of the 
diagram from Figure \ref{g0-2pt},
$c=\pm$ is the sign of deformation
of $K'$ into $K$,
and 
%enhanced
%eventually with the number
$\epsilon$
%(standing for
is the sum of the signs at all solitary nodes of $\R C$
(located in the region specified by the diagram)
%and the sign $c=\pm 1$
uniquely determines the rigid isotopy class of
a real algebraic curve.
Here the allowed values for $\epsilon$
are $\epsilon=\pm1$ in the cases 20-22, 28 and 30;
$\epsilon=0,-2$ in the case 26;
$\epsilon=\pm 1, -3$ in the case 27;
and $\epsilon=\pm 1,\pm 3$ in the case 29.
\end{coro}
We omit $\epsilon$ from $N^c_\epsilon$ in
the case when the diagram $N$ admits only
one value for $\epsilon$ (e.g. if there are
no solitary nodes at all). 
%Note that in the case $g=0$ we have $r=6$ and in the case $g=1$ 
%we have $r=5$ so that the hypothesis of Lemma \ref{rkprime} holds.
\ignore{
The sign $\pm1$ on the diagrams of Figure \ref{g0-2pt}
indicates the sign of the solitary real point of $\R C$
from $w(\R C)$.
The symbol $\epsilon$ stands for the sum of
the signs at all solitary nodes of $\R C$.

In Figure \ref{g0-2pt}.20-22, 28 and 30
we have $\epsilon=\pm 1$ while the solitary node
is unique and placed according to the picture.
In Figure  \ref{g0-2pt}.26
we have $\epsilon=0,-2$.
There
might be two solitary node in the same component
of $\rp^2\setminus\R C$ or no solitary nodes
(if $\epsilon=0$).
In Figure  \ref{g0-2pt}.27 we have $\epsilon=-3$
or $\epsilon=\pm 1$.
In Figures  \ref{g0-2pt}.29 
we have $\epsilon=\pm 3$ or $\epsilon=\pm 1$.
According to the value of $\epsilon$ there might be
one or three solitary nodes in the same component
of $\rp^2\setminus\R C$.
In particular, we claim that the rigid isotopy
type of the spatial curve $\R K$ depends
only on the number of the corresponding curve
$\R C$ in Figure \ref{g0-2pt} and the value of $\epsilon$
even if the number of solitary nodes of $\R C$ may vary.
%In Figure  \ref{g0-2pt}.29
%we may have one or three solitary real nodes
%in the corresponding curves while
%$\epsilon=\pm 1,\pm 3$.
%If $\epsilon$ is placed inside an oval then there is a single solitary node
%inside. If $\epsilon$ is outside there is one or several solitary nodes
%outside.
%In Figure \ref{g0-2pt}.20-22 $\epsilon=\pm 1$ and the solitary node
%is unique.
}
\begin{proof}
Once we ignore the sign of solitary nodes
(i.e. the $\sigma$-data in virtual diagrams
$(\R C;\R D^-,\tau,\sigma)$,
%in 
Figure \ref{g0-2pt} lists all triples
%is made of all possible virtual diagrams 
$(\R C;\R D^-,\tau)$ on
a generically immersed quartics $C$
with not more than 6 odd arcs
%of $\R D_+$ and two point of $\R D_-$
%after taking the quotient by
up to the equivalence generated
by all moves from Proposition \ref{prop-moves}.
We refer e.g. to \cite{DeMello} for
classification of generic real quartics $C$.

Suppose that $C$ does not have elliptic nodes.
In this case any 3D-explicit divisor
defines a nodal diagram unless $C$ has
a pair of conjugate nodes and $D$
is chosen so that the corresponding
lift $K'$ also have a pair of conjugate nodes.
This means that the real meromorphic function 
corresponding to $D-H_0$, where $H_0$ is
the divisor cut on $C\subset\p^2$
by the infinite axis of $\p^2$ has
real values at a fixed pair
of distinct and non-conjugate
points of $\CC C\setminus\R C$.
It is easy to see that the divisors
with this property form a codimension 2
subspace in the connected
real 8-dimensional space of all 3D-explicit
divisors corresponding to the same virtual
diagram. 

%If the planar rational curve $\R C$ has a pair
%of imaginary (complex conjugate)
%nodes
%%in the same component of $\rp^2\setminus\R C$
%then it can be degenerated to a solitary tacnode
%and further to a pair of
%complex conjugate imaginary
%solitary real nodes
%by the uniqueness of the rigid isotopy type
%for each homeomorphism class of the pair
%$(\rp^2,\R C)$.
%This deformation lifts to an isotopy in $\rp^3$
%once we choose a deformation of the divisor
%$\C D$ in a generic way.
%We claim that if the 
%solitary nodes have opposite signs then
%the deformation can be lifted to a deformation
%of $\R K'\subset\rp^3$ so that no singular points
%other than $p$ appear.

Our next claim is that if $(C,D)$
is a 3D-explicit nodal diagram such that $C$
has a pair of conjugate nodes and a
hyperbolic node then
there exists a path $(C_t,D_t)$, $t\in[0,1]$,
$(C_0,D_0)=(C,D)$ such that 
$(C_t,D_t)$ are 3D-explicit nodal diagrams,$C_1$
has three hyperbolic nodes,
$C_{\frac12}$ has a tacnode with two hyperbolic
branches while the curve in $\p^3$ corresponding
to $(C_{\frac12},D_{\frac12})$
by Proposition \ref{reconstruct-rkprime}
is smooth outside
of the projection point $p$.
The family $(C_t,D_t)$ determines
a rigid isotopy between
the corresponding algebraic knots.

To prove the claim
we note that there are two types of $(C,D)$
with such properties, and each can be obtained
by perturbation of two ellipses in $\p^2$ 
intersecting transversally at two
real points. 
%As we have already seen,
It is sufficient to prove the claim for one
representative curve in each type. 
The perturbations smooth one
of the real transversal intersection points
in two possible ways and keep the remaining
one real and two imaginary points of the intersection
of the ellipses.
Both perturbations can be included in a one-parametric
family of pairs of ellipses so that when $t$ increases
the ellipses become tangent and then intersect transversely in 4 distinct real points producing a family $C_t$
of rational nodal quadrics. We define $D_t$ so that 
%it stays disjoint from the nodes
%This ensures
%and that
$(C_t,D_t)$ is 3D-explicit,
and thus
the corresponding curve $K'_t\setminus\{p\}\subset\p^3$
is smooth.

Lemma \ref{enode} %and \ref{e2node}
allows us to reduce consideration
of nodal diagrams with elliptic nodes to those
without elliptic nodes and thus finishes the proof.
In particular, cases 26 with $\epsilon=+2$
as well as 27 with $\epsilon=+3$ can not
appear as the corresponding diagrams with
hyperbolic nodes have more than 6 odd arcs.
\end{proof}
 
\begin{lem}\label{enode}
Any real smooth rational sextic $K\subset\p^3$
is rigidly isotopic to a curve obtained
from a nodal diagram such that the underlying
real quartic rational curve does not have
elliptic nodes. 
%Suppose that $(C,D)$ is a 3D-explicit nodal diagram,
%where $C$ is a curve with an elliptic node.
%Then there exists a path $(C_t,D_t)$,
%$t\in[0,1]$,
%$(C_0,D_0)=(C,D)$ such that 
%$(C_t,D_t)$ is a 3D-explicit nodal diagram
%for $t\neq\frac12$, $C_{\frac12}$ is a cuspidal
%curve while $C_1$ has one elliptic node less than $C_0$.
%Furthermore, all curves $K'_t\subset\p^3$
%defined by $(C_t,D_t)$
%with the help of Proposition \ref{reconstruct-rkprime}
%have no singular points other than $p$.
\end{lem}
\begin{proof}
%We prove the lemma in two steps.
%In the first step we deform $(C,D)$ within
%the class of nodal diagrams so that all 6 points
%of $D$ become real.
%Assume that $C$ has an elliptic node.
%Suppose that $D^-\subset\R C$.  
%First we note that we may assume that
%$D^-\cap\R K=\emptyset$.
Let $q\in \R C$ be an elliptic node in
the nodal diagram $(C,D)$ of $K$.
%The two other nodes may be either both
%real or form a complex conjugate pair.
%In the latter
%
The image of $C$ under the quadratic transformation
centered in the nodes of $C$ is a conic intersecting
the coordinate line corresponding to $q$ in 
two imaginary points. A deformation of this conic
to a conic intersecting this line in two real
points give a deformation $C_t$, $t\in[0,1]$,
of $C=C_0$ into a rational
nodal quartic $C_1$ that changes $q$ into a hyperbolic
node, and leaves the other nodes of $C$ unchanged. 
This deformation can be extended to a deformation
$D_t$ of $D=D_0$ so that $D_t$ remains
disjoint from the nodes of $C_t$ and neither $D^+_t$
nor $D^-_t$ has multiple points.

This gives a deformation $K'_t\subset\p^3$
of sextic curves. If $K'_t\setminus\{p\}$
is nonsingular then the smooth curves obtained by
coherent deformations of $K'_0$ and $K'_1$
are rigidly isotopic while the number
of elliptic points of $C_1$ is less than that of $C_0$,
so that we may proceed inductively.
Suppose that $K'_t\setminus\{p\}$ is singular
for $t=\epsilon$ and nonsingular for $t\in [0,\epsilon)$.

Note that the singularity of $K'_t\setminus\{p\}$
must sit
over an elliptic node $s$ of $C_t$ since $D$ is
3D-explicit. Exchanging the roles of $s$ and $p$
if needed,
we may assume that $p$ was elliptic, i.e.
$D^-\cap\R K=\emptyset$.
In such case $\R K\setminus\nu^{-1}(\Sigma)$
consists of not more that 6 arcs,
and $D^+$ maybe deformed in a family $D_t$
of divisors in the same curve $K$
so that $D^-_t=D^-$, $D^+_t\cap D^-_t=\emptyset$,
$D^+_t\cap\nu^{-1}_1(\Sigma)=\emptyset$, $t=[0,1)$,
while $D^+_1=\nu^{-1}(\Sigma)$, so that
$K'_1\subset\p^3$ is a quadrinodal sextic curve,
i.e. a curve with 4 distinct nodes.
Note that
%the Bezout theorem ensures that
by our construction these 4 nodes
are not coplanar and the curve $K'_t$ is
not contained in any plane of $\p^2$.
If there are nodes forming a complex
conjugate pair then we can proceed
as in the proof of Corollary \ref{coro-2ptg0}
deforming this pair to a pair
of hyperbolic nodes.

If all nodes are real
(elliptic or hyperbolic) then
we choose the coordinates in $\p^3$
so that the intersections of the coordinate
hyperplanes correspond
to the nodes of $K'_1$, 
the cubic transformation 
\begin{equation}\label{cucr}
(x_0:x_1:x_2:x_3)\mapsto (\frac 1{x_0}:\frac1{x_1}:\frac 1{x_2}:\frac 1{x_3})
\end{equation}
maps $K'_t$ to a conic in $\p^3$. Thus all
quadrinodal curves corresponding to the same
planar nodal quartic $C$ are isotopic.
But all elliptic nodes of $C$ can
be simultaneously deformed to hyperbolic nodes
as we can see
through consideration of conics in $\p^2$
tangent to coordinate lines.%hyperbolic 

It remains to prove that
we can deform $(C,D)=(C_0,D_0)$ to $(C_1,D_1)$ so that all
$(C_t,D_t)$, $t\in[0,1)$ are nodal diagrams while $K'_1\subset\p^3$ is quadrinodal.
We start with a deformation of $D$ on the same curve $K$ as considered above.
If there exists $\epsilon<1$ with singular $K'_\epsilon\setminus\{p\}$ then its singularity
must be an elliptic node $e$. Consider a plane section
$H$ of $K'_\epsilon$ disjoint from $\{p\}$ and such
it passes through $e$ and such that the corresponding
divisor is 3D-explicit (we can do that since 
there are not more than 2 hyperbolic nodes of $C$).
Deforming the plane section from $D^+_\epsilon$
to $H$ gives us a family of intermediate nodal diagrams
$(C_t,D_t)$, $t\in [\epsilon,\epsilon_H)$,
$\epsilon<\epsilon_H<1$ while $D_{\epsilon_H}^+$ contains
$\nu_{\epsilon_H}^{-1}(e)$. We continue the process
by deformation of the effective divisor
$D_{\epsilon_H}^+-\nu_{\epsilon_H}^{-1}(\Sigma_{\epsilon_H}\setminus\{e\})$ until we arrive to $D_1=\nu^{-1}_1(\Sigma_1)$. Finally we perturb the family $(C_t,D_t)$
slightly to ensure that it consists if nodal diagrams
for $t<1$.

\ignore{
\newcommand{\ch}{\operatorname{Ch}}
Let us recall that a rational nodal quartic
$C\subset\p^2$ is determined (up to
a multiplicative translation) by 
its {\em chord diagram} on $\p^1$ defines as three
disjoint pairs of points $\ch_j\subset \p^1$, $j=0,1,2$.
Indeed, three pairwise unions of $\ch_j$ define
three divisors of degree 4 on $\p^1$,
%with empty common intersection
and thus a rational curve of degree 4
in $\p^2$ with nodes at the three points
of intersection of coordinate axes.
An elliptic node corresponds to a pair
of points interchanged by $\conj$. 

The divisor $D^+$ is given by 6 $\conj$-invariant
points on $\p^1$. Clearly,
a pair of conjugate points of $D^+$
may be deformed to a pair of real points that
are close to each other and disjoint
from $D^-$ and $\ch_j\cap\rp^1$.
We claim that we may deform $D^+$ and $\ch_0$
simultaneously so that ..
simultaneously 
}
\end{proof}

\begin{figure}[h]
\includegraphics[width=125mm]{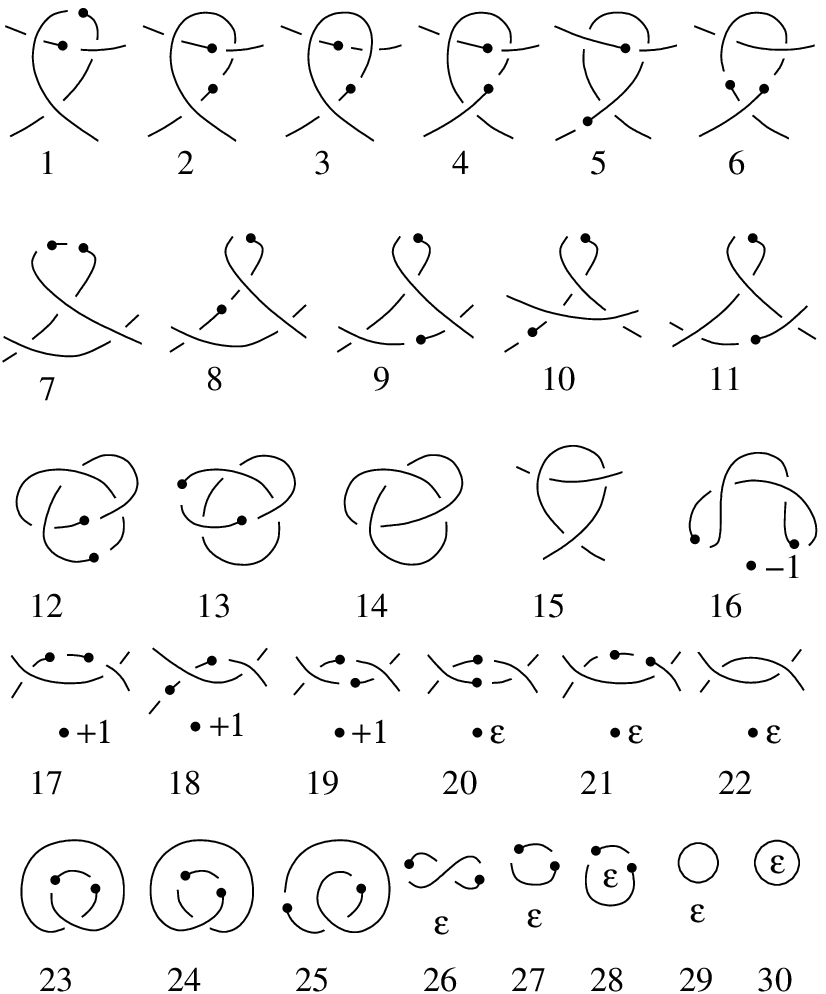}
\caption{
Equivalence classes of the diagrams of $\R K'$ with respect to the moves
of Proposition \ref{prop-moves}.
\label{g0-2pt}}
\end{figure}

\subsection{Quadrinodal spatial rational sextics
and isotopy}
The quadrinodal curves we considered in the previous
subsection are also useful for 
finding isotopies among the curves
listed in Corollary \ref{coro-2ptg0}.
%Our main source for such isotopies is {\em real quadrinodal rational sextic curves},
%i.e. rational real algebraic
%curves $\R J\subset\rp^2$ of degree 6 with four real double points
%(with distinct tangent directions
%of the branches).
%We refer to each double point of $\R J$ as the {\em node}.
%We call a real node {\em hyperbolic}
%if it corresponds to 
%an intersection point of two real branches and {\em elliptic}
%otherwise (if it is a solitary real point).
Let $J\subset\p^3$ be a quadrinodal rational
sextic with 
the nodes in $(1:0:0:0)$, $(0:1:0:0)$, $(0:0:1:0)$
and $p=(0:0:0:1)$.
%the projection point $p$ is one of its nodes.
\ignore{
If we project $\R J$ to $\rp^2$
from one of its double point the result is a quartic curve $\R C\subset\rp^2$  enhanced
with a real divisor $\C D=\C D_+-\C D_-$ on its complexification. 
As before, $\C C_+$ and $\C D_-$ are effective divisors with $\deg\C D_+=6$
and  $\deg\C D_-=2$. Here $\C D_-$ is formed by the tangent lines
to the two branches of $\C J$ at the double point chosen as the center $p$ of projection
while $\C D_+$ is given by a choice of the reference ``coordinate plane" not
passing through $p$.

It is convenient to choose the four homogeneous coordinate planes in $\rp^3$
to be the four planes passing through the four triples
formed from the four
double points of $\R J$.
By the Bezout theorem there are no
intersection points of $\R J$ with these coordinate planes
other than the double points.
We call such a coordinate system
{\em compatible with the nodes of $\R J$.}

Recall that for each of the four nodes of $\R J$
we may consider its resolution into a positive or negative
crossing as in Figure \ref{perturb}. 
%Furthermore, we may choose to leave the node unresolved.
%We refer to it 
}
\begin{prop}
For every choice of signs
for some nodes of a
real rational quadrinodal sextic curve $J$
there exists a deformation of $\R J$ in the class or real rational
sextics
resolving those nodes according to the chosen signs
and keeping all the other nodes unperturbed.
\end{prop}
\begin{proof}
%Choose a node $p\in\R J$ as the center of the projection $\pi_p:\R J\setminus\{p\}\to\rp^2$.
%Consider the presentation given by Proposition \ref{cd-rkprime}.
%Our choice of the coordinate plane $\R H$ is such that the six points
%of the divisor $\C D_+$ are located at the nodes of $\R C$%.
The curve $J$ corresponds to a nodal rational planar
curve $C$ and the divisor $D^+$
consisting of 6 points in the normalization $K$ of $C$
corresponding to 3 nodes of $C$
so that each node corresponds to a pair of points
in $D^+$. 
%We have two points of $\C D_+$ at each node, one on each branch.
Moving one of the points in the pair in
an appropriate direction we get
the deformation of $\R J$ with the chosen sign.
%is determined by 
%the condition 
(In the case when the pair consists of complex
conjugate points we move the second point 
in a complex conjugate way.)
If our choice of signs keeps some of the nodes
unresolved then without loss of generality we
may assume that we do not resolve $p$.
If our choice resolve all the nodes
%This procedure resolves all points but $p$. If the choice of signs
%is made at all four points (including $p$)
then we apply Lemma \ref{rkprime}
to resolve $p$.
\end{proof}

The normalization of $\R J$ is topologically a circle.
Hyperbolic nodes of $J$ may be encoded on this circle
by means of the so-called
{\em chord diagram}: we draw a chord
%on the circle $S=\rp^1$
connecting
each pair of points that gets identified to a node of $\R J$.
Real solitary nodes of $\C J$ are ignored
in the chord diagram $S$,
their number is equal 
to 4 minus the number of chords.
As the space of $n$ distinct
pairs of complex conjugate points in $\CC\p^1$
is connected we obtain the following statement.
\begin{prop}
The space of real quadrinodal rational sextic curves in $\p^3$ (considered up to projective linear transformations) corresponding to the same combinatorial diagram of real chords
is connected.
\end{prop}
Here we consider only quadrinodal curves with all
nodes real, i.e. without complex conjugate pairs
of nodes.
\subsection{Completing the classification}
\begin{proof}[Proof of Theorems \ref{thm-g0t} and \ref{thm-g0r}]
By Corollary \ref{coro-2ptg0} to deduce the classification we need
to identify the curves obtained from the diagrams of Figure \ref{g0-2pt}
that are rigidly isotopic.
For this we use eleven real rational quadrinodal curves
given by the chord diagrams A through K depicted on
Figure \ref{4dp}
%the chord diagrams A through K depicted on Figure \ref{4dp}
as well as the 3-nodal curve L from Figure \ref{4dp-l}.

\begin{figure}[h]
\includegraphics[width=78mm]{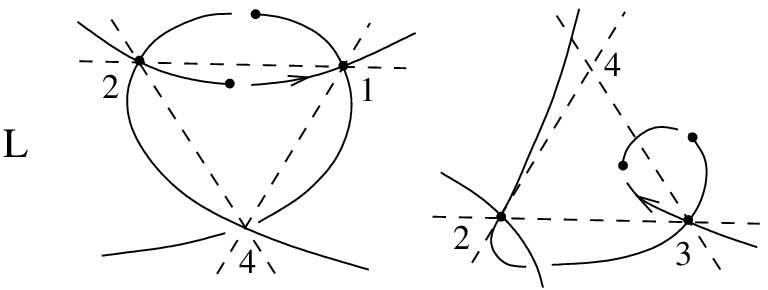}
\caption{
A trinodal rational sextic curve.
\label{4dp-l}}
\end{figure}
\begin{figure}[h]
\includegraphics[width=90mm]{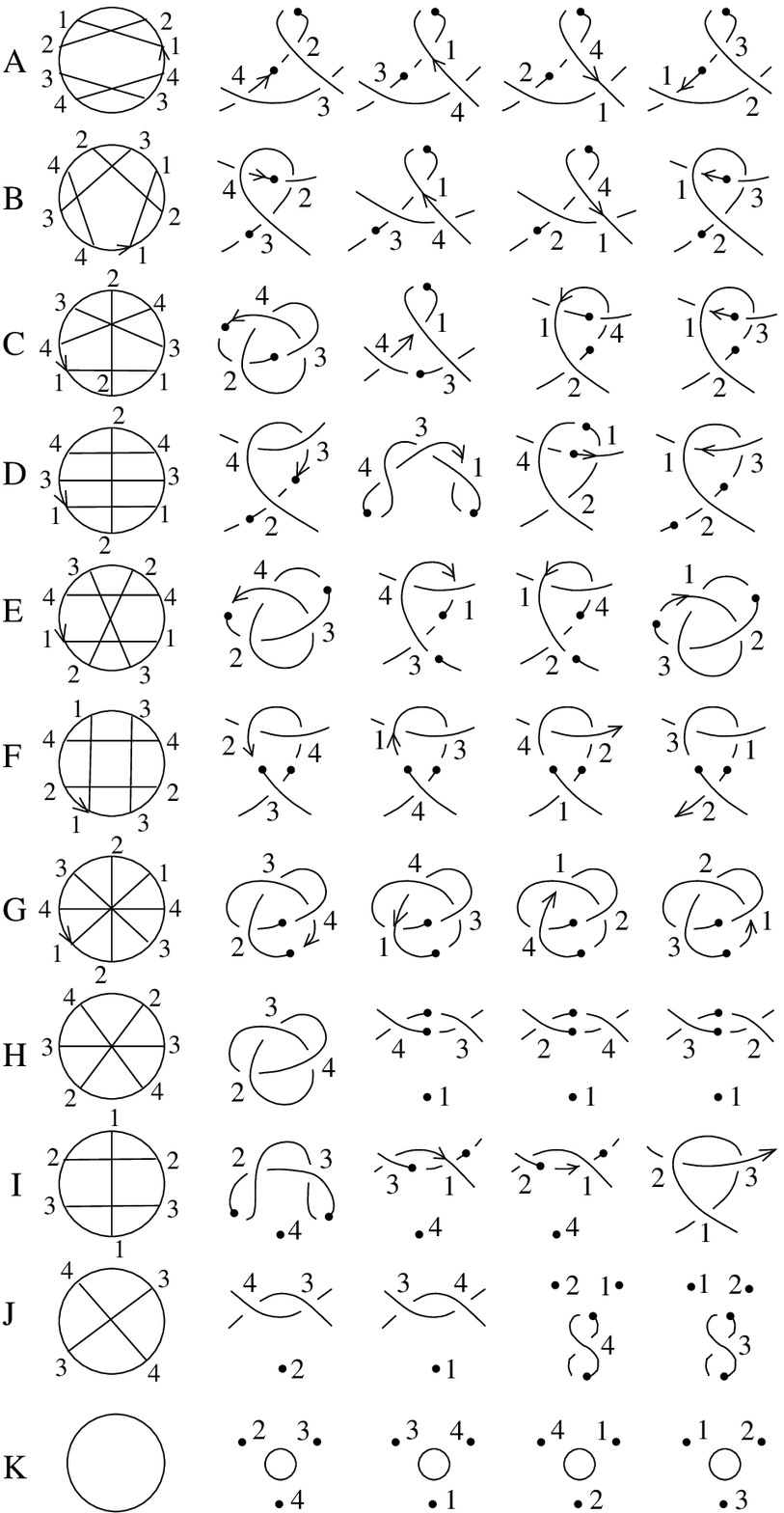}
\caption{
Eleven quadrinodal sextic curves.
\label{4dp}}
\end{figure}
The quadrinodal curves are depicted along with 4 projections
from each of its nodes (numbered 1 through 4).
For convenience we indicate the diagram
of the one-nodal curve $\R K'$
obtained by {\em positive} resolution of all nodes except for the
projection point.
To identify different projections of the same quadrinodal curve we indicate
the same arc connecting two of the nodes.

Also we depict the lower and upper half arcs near the points of $\C D_-$
for the corresponding diagrams. Note that this choice is determined
(up to simultaneous reversal) by the following rule.
Let us connect the two points of $\C D_-$ by an arc in
the curve $\R C\subset\rp^2$ and compute the number of points of $\C D_+$ contained
on this arc. If the arc crosses the infinite line of the projective plane of
the diagram odd number of times
then we add one to this number. If the result is odd
then the arc connects a lower half-arc to an upper half arc. Otherwise it connects
the half-arcs of the same kind.

%To construct the trinodal curve ``L'' we use Proposition \ref{cd-rkprime}
%and place four of the points of $\C D_+$ to the normalizations of the nodes 1 and 2.
%The remaining two points are placed in the imaginary
%domain $\C C\setminus\R C$.
%Note that this diagram cannot be obtained by a resolution of a quadrinodal curve
%as the half arcs distribution near $\C D_-$ violate the parity condition
%considered above.
The tri-nodal curve in Figure 10 (curve ``L'') is parameterized by
$$
  t\mapsto(
  p_1 p_3 p_5 p_6 p_7 p_8 :
  p_1 p_2 p_3 p_4 p_5 p_7 :
  p_2 p_4 p_6^3 p_8 :
%                 ^
  p_1 p_2 p_4 p_5 p_6 p_8 ),
$$
where $p_i=t-t_i$ and
$t_1 < t_2 < \dots < t_8$.
\begin{figure}[h]
\includegraphics[width=125mm]{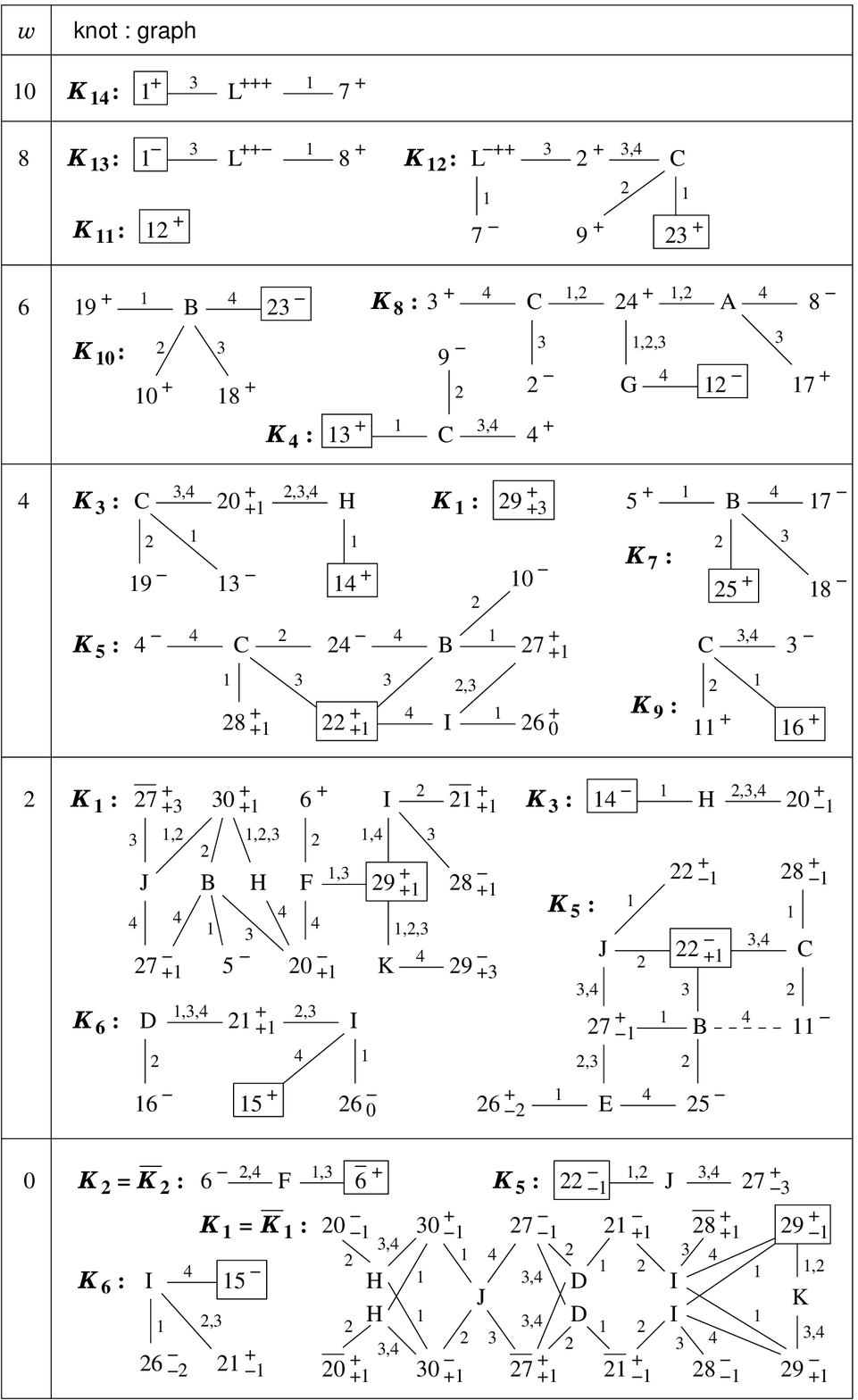}
\caption{
Graphs of rigid isotopy equivalence.
\label{graphs}}
\end{figure}

Figure \ref{graphs} indicates which of the
curves obtained from the diagrams of Figure \ref{4dp}
%by positive or negative resolution of the node
are rigidly isotopic.
Namely, Figure \ref{graphs} lists bipartite graphs with two type
of vertices: numeric and alphabetic.
The number of the numeric vertices refers to the diagram number
from Figure \ref{g0-2pt} as well as the sign used for the perturbation of the node.
Each such diagram encodes the equivalence class with respect to the moves
of Proposition \ref{prop-moves}.

The letter of the alphabetic vertices refers to the multinodal curves from Figure \ref{4dp}.
Each edge is labeled by
the number of the node of the corresponding
multinodal curve that becomes the projection point in the diagram for the adjacent numeric
vertex.
The signs of the resolution used in the graph are determined by the signs
at the adjacent numerical vertices.

We see that each pair consisting of a curve from Figure \ref{deg6c} and 
the non-negative value of its Viro invariant $w$ corresponds to a connected subgraph
and that each resolution of a curve from Figure \ref{g0-2pt} is contained
in one of the subgraphs. Different knots from Figure \ref{deg6c} are topologically
different as knots in $\rp^3$, see \cite{JuViro}.
Figure \ref{ak0} provides topological identification of the boxed resolved diagram
in Figure \ref{graphs} and the corresponding knot type from Figure \ref{deg6c}.
\begin{figure}[h]
\includegraphics[width=115mm]{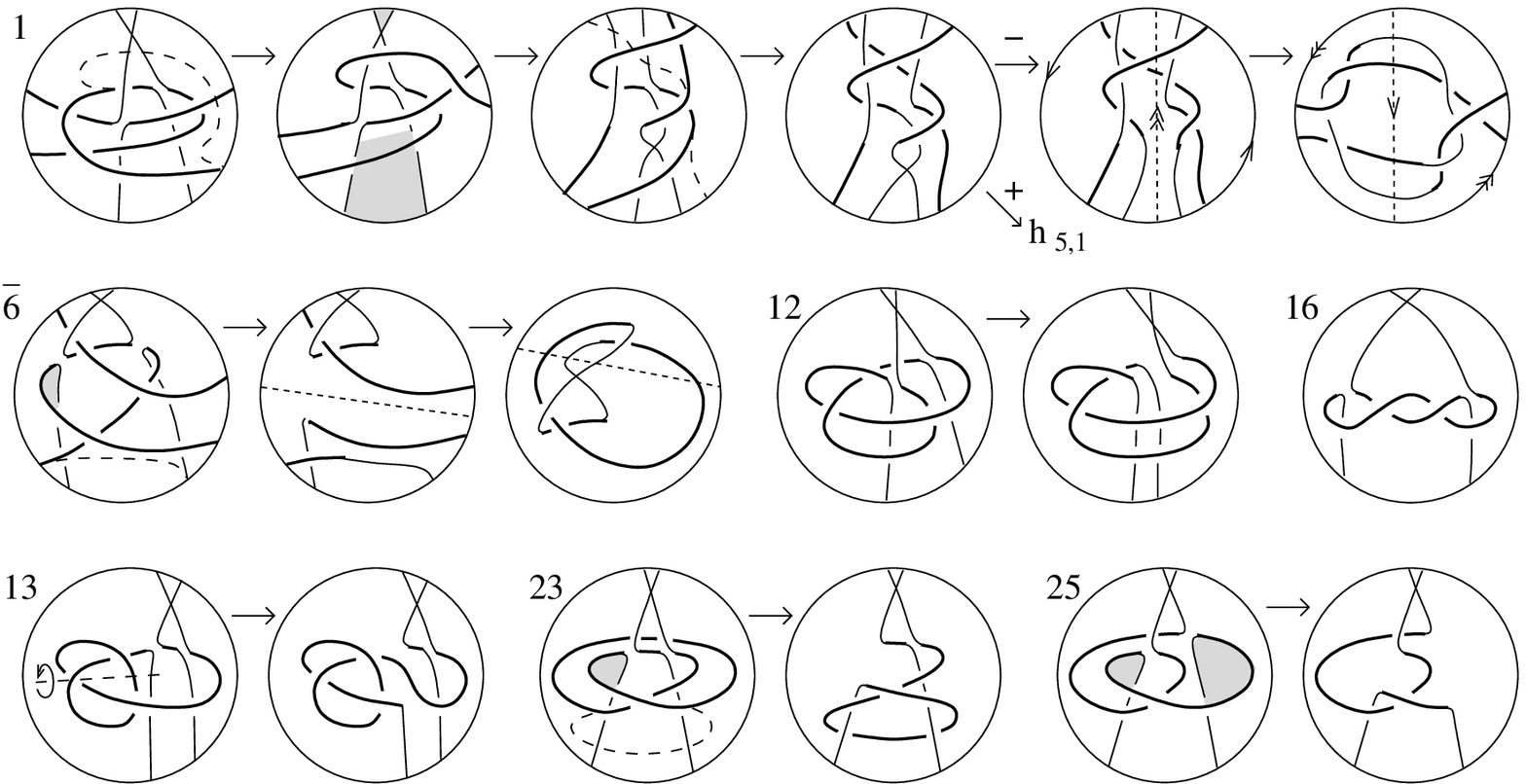}
\caption{
Identifications between the knots from Figure \ref{deg6c}
and the boxed diagrams of Figure \ref{graphs}.
\label{ak0}}
\end{figure}

A reflection (an orientation-reversing automorphism) in $\rp^3$ reverses the knot together
with its Viro invariant.
Thus each curve with positive $w$ corresponds to two distinct rigid isotopy types
under reflection
while a curves with $w=0$ may correspond to one or two types.
We have four types of knots with $w=0$: $K_1$, $K_2$, $K_5$ and $K_6$.
The knots $K_5$ and $K_6$ are chiral: they are not topologically isotopic
to their reflection in $\rp^3$ as the two components of their inverse images
under the universal covering ${\mathbb S}^3\to\rp^3$ have non-zero linking number.
In the same type the knot types $K_1$ and $K_2$ are amphichiral.
Furthermore, the corresponding real algebraic sextics are rigidly isotopic
as their reflections appear in the same components of the graphs from Figure \ref{graphs}.
Altogether we get 38 rigid isotopy types in Theorem  \ref{thm-g0r}.
We get 14 topological types of Theorem \ref{thm-g0t}
by taking out the Viro invariant information from the data. 
\end{proof}

%% file: d6g1.tex
\section{Elliptic knots and links of degree 6}
\begin{figure}[h]
\includegraphics[width=100mm]{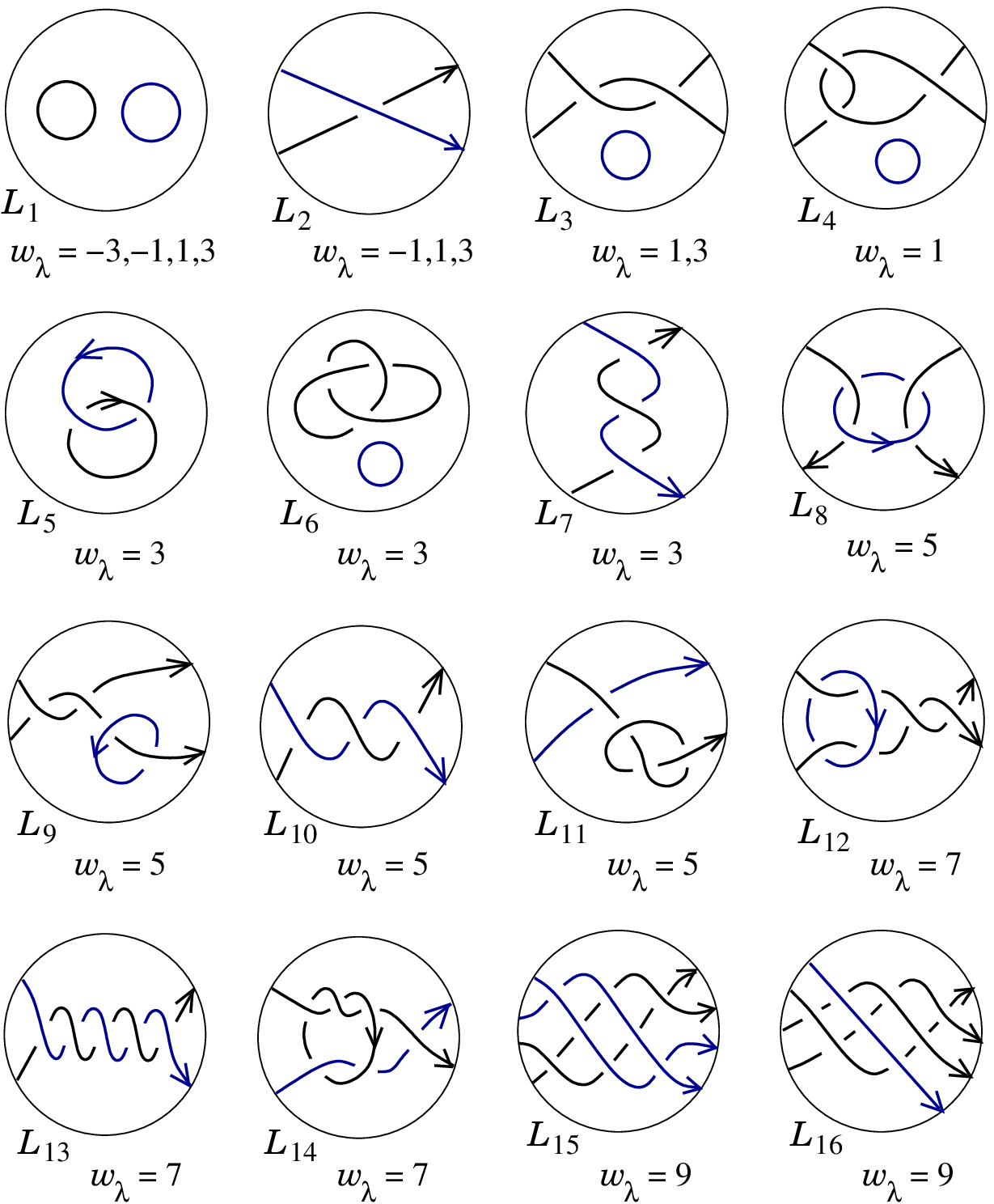}
\caption{
Two-component real algebraic links of degree 6 and genus 1.
\label{deg6l}}
\end{figure}

%Finally we treat the case when $L\subset\p^3$
%is an {\em elliptic} (genus 1) smooth real curve 
%of degree
\newcommand{\ch}{\operatorname{ch}}
%\subsection{Statement of classifications}
\begin{thm}\label{thm-g1t}
There are 16 topological isotopy types (homeomorphism classes
of pairs $(\rp^3,\R L)$)
of elliptic (genus 1) smooth real algebraic curves $K$ of
degree 6 %embedded
in the projective space $\rp^3$ in the case
when $\R K$ is a two-component link.
% considered up to
%diffeomorphism of the pair $(\rp^3,\R K)$.
Namely, Figure \ref{deg6l} lists the links.

There are 4 such types in the case when $\R L$ is
connected. Namely, the types $K_1,K_3,K_5,K_6$
from Figure \ref{deg6c} are realizable by
smooth elliptic sextic curves in $\p^3$.
\end{thm}

\begin{thm}\label{thm-g1r}
There are 40 rigid isotopy types of rational real algebraic curves of degree 6
embedded in the projective space $\rp^3$ in the case when $\R K$ is a two-component
link.
Namely, each curve depicted on Figure \ref{deg6l} enhanced with 
a choice of the listed value for $\wl$ gives rise to
one rigid isotopy class of $\R K$ in the depicted link type. 
Furthermore, simultaneous reflection of $\R K$ in $\rp^3$ and changing
the sign of $\wl$ gives a new rigid isotopy type with the exception
of when $\R K$ is a 
topologically trivial link (the first case
in Figure \ref{deg6l}).

There are 12 such types in the case when $\R K$ is
connected. Namely, we may have $\wl=\pm 3,\pm 1$
for the types $K_1$; $\wl=3$ for
$K_3$; $\wl=1,3$ for $K_5$; and $\wl=1$ for $K_6$
(see Figure \ref{deg6c}).
\end{thm}

\subsection{Binodal planar elliptic quartics}
%and trinodal spatial elliptic sextics}
To prove Theorems \ref{thm-g1t} and \ref{thm-g1r}
we use a technique similar to that developed
in the previous section with some modifications.
Lemma \ref{rkprime} and Propositions 
\ref{pirkprime} and \ref{reconstruct-rkprime}
reduce the proof to studying of nodal diagrams
$(C,D)$, where $C$ is a nodal elliptic curve
 of degree 4 while $D=D^+-D^-$ is linearly
equivalent to a hyperplane section of $C$.
%$\deg D^+=6$, $\deg D_-=2$.
Namely, the nodal diagram $(C,D)$ defines
a spatial curve $K'\subset \p^3$
with the node $p=(0:0:0:1)\in \p^3$
that can be resolved in two ways to a smooth
spatial elliptic sextic. Once again, we
can distinguish between positive and negative
resolution of $p$. 
If $p$ is the self-intersection
of the same component of $\R K$
then we may use any of its orientations
to determine the sign.
If $p$ is the intersection of different
components then $K$ is necessarily of type I
and we may use a complex orientation of 
$\R K$ to determine the sign. 

The singular set $\Sigma\subset C$ consists
of two nodes. Thus the four points
of $\nu^{-1}(\Sigma)\subset K$ 
give a hyperplane section of $C$  
through the normalization map $\nu:K\to C$.

Let $E$ be an abstract elliptic curve
(not embedded to any projective space)
enhanced with an antiholomorphic involution
$\conj$ with non-empty fixed locus $\R E$.
Let $\ch_j\subset E$, $j=1,2$, be two disjoint pairs
of points with $\conj(\ch_j)=\ch_j$.

\begin{lem}\label{2chords}
If the divisors formed
by $\ch_1$ and $\ch_2$ are not linearly 
equivalent in $E$ then there exists a
planar nodal
quartic curve $C\subset\p^2$
with two nodes $q_j\in C$, $j=1,2$,
with $(K;\nu^{-1}(q_1),\nu^{-1}(q_2))$
is isomorphic to $(E;\ch_1,\ch_2)$.

Vice versa, for any nodal irreducible quartic
curve $C$ with two nodes $q_j\in C$, $j=1,2$,
the divisors $\nu^{-1}(q_1)$
and $\nu^{-1}(q_2)$ are not linearly equivalent.
\end{lem}
\begin{proof}
Consider the projective linear system
%of divisors in the equivalence class of
$|\ch_1+\ch_2|$ and a point
$r\in E\setminus (\ch_1\cup\ch_2)$.
As $\dim |\ch_1+\ch_2|=3$
%his system is of projective rank 3. Choose
there exist a unique point
%$r\in E\setminus (\ch_1\cup\ch_2)$.
%The two points
$r_j$, $j=1,2$, such that
$r_j+r+\ch_j$ is equivalent to $\ch_1+\ch_2$.
Since $|\ch_1|\neq|\ch_2|$ we have $r_1\neq r_2$,
and $\ch_1+\ch_2, r_1+r+\ch_1, r_2+r+\ch_2$
generate a 2-dimensional subsystem in 
$|\ch_1+\ch_2|$ and a map of
$E$ to $\p^2$. Note that
the image of $E$ has two nodes corresponding to
$\ch_j$. Thus
this subsystem
contains $s_j+s+\ch_j$ for with some $s_j$
for every
$s\in E\setminus (\ch_1\cup\ch_2)$. 
Therefore, the resulting linear system
is independent of the choice of $r$.

\ignore{
the 2-dimensional subsystem generated by 

Choose a linear system of degree 3 to present $E$
as a planar cubic. Choose homogeneous
coordinates in $\p^2$
so that the pairs $\ch_j$ are contained in the
coordinate lines $\{z_j=0\}$. Perturbing our choice
of linear system if needed we may assume that
the coordinate lines $\{z_j=0\}$ intersect
$E$ in three distinct points. Let
$r_j\in \{z_j=0\}\cap E$ be the remaining
intersection point.
Choose the third coordinate line to contain
$r_1$ and $r_2$ and make the quadratic transformation
according to this coordinate choice.
}

For the converse it suffices to take
a generic point in a nodal
elliptic quartic curve
$C\subset\p^2$ and connect it with the nodes $q_j$
of $C$ by lines. The lines will intersect $C$
at distinct fourth points. Distinct single points
are not linear equivalent since our normalization $K$
is not rational.
\end{proof}

This allows us to work with planar nodal elliptic quartics
almost as freely as with planar nodal rational quartics.

\begin{coro}\label{binod}
There are
10 distinct rigid isotopy classes of 
real elliptic quartic curves in $\p^2$
%with two hyperbolic nodes
in the case when the normalization of the real
locus consists of two components:
5 with both nodes
hyperbolic, 2 with one hyperbolic and one elliptic
node, 1 with two elliptic nodes and 2 with a
pair of complex conjugate nodes.
There are 5 distinct classes in the case when the
normalization is connected: 2 with both nodes
hyperbolic, 1 with with one hyperbolic and one elliptic
node, 1 with two elliptic nodes and 1 with a
pair of complex conjugate nodes.
\end{coro}
Here by rigid isotopy we mean a deformation 
in the class of (irreducible) real elliptic binodal
quartics in $\p^2$.
\begin{proof}
It is convenient to think of $\ch_j$ as a {\em chord}
connecting two points of the real locus $\R E$.
We have three deformation classes if each chord
connects two points from the same component of $\R E$,
and
a single class if one or both chords connect
different components of $\R E$.
Each class is unique to automorphism and deformation
of $(E;\ch_1,\ch_2)$ in the class of triples with
$|\ch_1|\neq|\ch_2|$.
If $\R E$ is disconnected then complex conjugate
nodes may correspond to intersections of the same
or different components, this gives us two classes.
\end{proof}

\subsection{Nodal diagrams %and their moves
in the
elliptic case}
\begin{figure}[h]
\includegraphics[width=80mm]{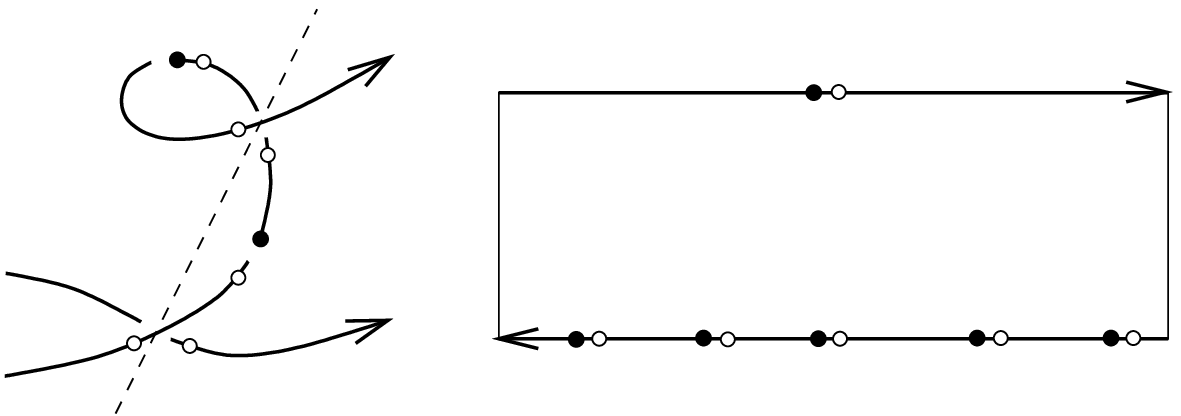}
\caption{
Six odd arcs and the Abel theorem.
\label{lemdd}}
\end{figure}

Propositions \ref{prop-diagram} has the following
counterpart for the case of genus 1.
\begin{prop}\label{prop-diagram-g1}
Suppose that $(\R C;\R D^-,\tau,\sigma)$
is a virtual nodal diagram such that
all nodes of a nodal elliptic quartic curve $\R C$
are hyperbolic.
The diagram $(\R C;\R D^-,\tau,\sigma)$
is not realizable by a 3D-explicit nodal diagram $(C,D)$
%of
%a rational real curve $\R K'\subset\rp^3$ of degree 6
%with a double point at $p\in\rp^3$
%if the number of its odd arcs
%is at most 4. 
if the six-point set $D^-\cup\nu^{-1}(\Sigma)$
is contained in a connected component of $\R K$.
It is also not realizable if
$D^-\cup\nu^{-1}(\Sigma)\subset\R K$ and 
$\R C$ has a node corresponding to intersection of
distinct components of $\R K$.

In all other cases
$(\R C;\R D^-,\tau,\sigma)$
is realizable by a 3D-explicit nodal diagram $(C,D)$.
%If this number is more than 6 it is not realizable.
%
%If the number of  equal to 6 ..
%In the realizable cases a nodal diagram $(C,D)$
%for
%such a realizable $(\C C;\{x'',y''\})$
%may be chosen to be 3D-explicit.
Furthermore,
%the diagram $(C,D)$ may be chosen to be 3D-explicit,
%and
the virtual diagram $(\R C;\R D^-,\tau,\sigma)$
together with the sign $c$
of deformation of $K'$
determine the embedded real algebraic curve $K\subset\p^3$ up to rigid isotopy.
%any two real rational curves of degree 6 
%obtained from curves corresponding
%to the same virtual nodal
%diagram with with nodal diagrams determining a 3D-lift
%corresponding
%to $(\C C;D^-,\tau,\sigma)$ with the same value of $c$
%the same virtual nodal diagram are 
%are rigidly isotopic.
%isotopic in the class of rational real curves with
%nodal diagrams determining a 3D-lift corresponding
%to the same virtual nodal diagram.
\end{prop} 
\begin{proof}
Suppose $(C,D)$ is a 3D-explicit nodal diagram.
%The
%number of odd arcs can not be greater than 6
%as $\deg D^+=6$.
The number of arc-components of
$\R K\setminus(\nu^{-1}(\Sigma)\cup D^-)$
is 6 if $D^-\subset\R K$, otherwise it is 4.
Suppose that we have 6 odd arcs.
 Then
all 6 points of $D^+$ are real, and each point
sits on its own arc of $\R K\setminus(\nu^{-1}(\Sigma)\cup D^-)$. If all these 6 points belong to
the same component of $\R K$ we have a
contradiction to the Abel theorem
as $|D^+|=|D^-+\nu^{-1}(\Sigma)|$, but
a divisor cannot be principal if it
can be presented as the boundary
of a collection of disjoint coherently
oriented intervals contained
in the same component of $\R K$.

If the images of different components of $\R K$
under $\nu$ intersect then by Corollary \ref{binod}
there are just two possibilities for $\R C$
up to rigid isotopy. In this case $K$ is of type I.
%so we choose a complex orientation on $\R K$.
Note that if $f:K\to\p^1$ is a real meromorphic
function with $f^{-1}(\rp^1)=\R K$ then $f$
has no real critical values, and a complex
orientation of $\R K$ can be obtained 
as the pullback of an orientation of $\rp^1$.
Applying this remark to the function defined
by $D^+-(D^-+\nu^{-1}(\Sigma))$
for these two possibilities we get a contradiction, see Figure \ref{lemdd}.
%If this number is
%equal to 6 then all points of $D^+$ are real
%and the
In the remaining cases existence and
uniqueness of a nodal diagram $(C,D)$
up to deformation follows from Lemma \ref{2chords}, where the black points represent $f^{-1}(\infty)$
(i.e. the divisor $D^- + \nu^{-1}(\Sigma)$) and the white points
represent $f^{-1}(\epsilon)$
with $\epsilon\gg 0$. 
\end{proof}

\begin{lem}\label{enode-g1}
Any real smooth rational sextic $K\subset\p^3$
is rigidly isotopic to a curve obtained
from a nodal diagram such that the underlying
real quartic rational curve does not have
elliptic nodes. 
\end{lem}
\begin{proof}
Suppose that $(C,D)$ is a 3D-explicit
nodal diagram, and $C$ has an elliptic node.
%Note that in such case
%$\R K\setminus(\nu^{-1}(\Sigma)\cup D^-)$
%consists at most of four arcs, so that we may
%choose $D^+$ so that $\deg (D^+\cap\R K)\le 4$.
Use Lemma \ref{2chords} to deform the corresponding
chord %of the normalization
to a tangent line, and then further to real chord.
Then the elliptic node gets deformed to a cusp,
and further to a hyperbolic node.
Deform $D$ in its linear equivalence class to extend
the deformation of $C$ to a 3D-explicit deformation
$(C_t,D_t)$
of $(C,D)$. We can do that since we may assume
that $\deg (D^+\cap\R K)\le 4$
since $\R K\setminus(\nu^{-1}(\Sigma)\cup D^-)$
contains at most of four arcs.

If during the deformation the curve
$K'_t\setminus\{p\}\subset\p^3$ has a singular
point then it must be an elliptic node. 
Thus exchanging the roles of $p$ and this node
we may assume that $p$ is elliptic, i.e. 
$D^-\cap \R K=\emptyset$. In this case
we may assume that $\deg (D^+\cap\R K)\le 2$
and keep two of the points of $D^+_t$
at the inverse image of the elliptic node of $C_t$
under $\nu$.
If $C_t$ has two elliptic nodes then $D^+_t$
can be chosen without real points and we may 
keep four points of $D_t^+$ at $\nu^{-1}(\Sigma)$
thus ensuring (after a perturbation)
that $K'_t\setminus\{p\}$ is smooth.
\end{proof}

\begin{figure}[h]
\includegraphics[width=90mm]{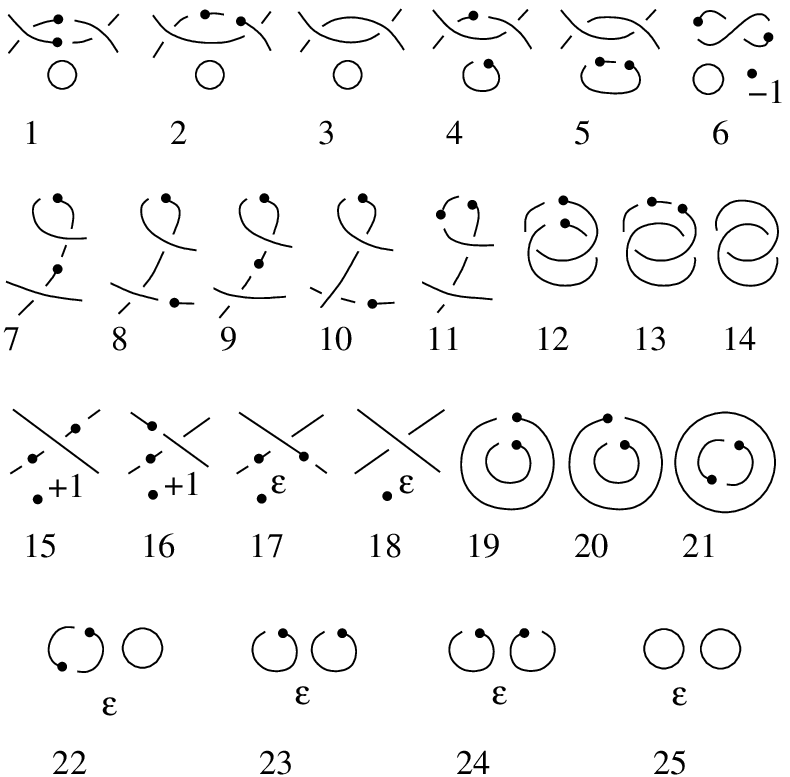}
\caption{
Equivalence classes of the diagrams of $\R K'$ with respect to the moves
of Proposition \ref{prop-moves}.
\label{g1-2pt}}
\end{figure}

\begin{coro}
\label{coro-2ptg1}
Any real elliptic curve $K$ of degree 6 in $\p^3$
is rigidly isotopic
to a curve obtained from one of the curves
whose virtual nodal diagrams
are listed on Figure \ref{g1-2pt}
by resolving them according to $c=\pm 1$ 
as in Lemma \ref{rkprime} if $K$ is of type I.

Each diagram $N^c_\epsilon$,
where $N$ is the number of the 
diagram from Figure \ref{g1-2pt},
$c=\pm$ is the sign of deformation
of $K'$ into $K$,
and 
%enhanced
%eventually with the number
$\epsilon$
%(standing for
is the sum of the signs at all solitary nodes of $\R C$
(located in the region specified by the diagram)
%and the sign $c=\pm 1$
uniquely determines the rigid isotopy class of
a real algebraic curve.
Here the allowed values for $\epsilon$
are $\epsilon=\pm1$ in the cases 17 and 18;
$\epsilon=0,-2$ in the cases 22 and 23;
and $\epsilon=0,\pm2$ in the cases 24 and 25.

If $K$ is of type II then it 
is rigidly isotopic to a curve obtained
from is rigidly isotopic
to a curve obtained from one of the curves
whose virtual nodal diagrams
are listed on Figure \ref{g0-2pt},
cases 20, 21, 22 (without further elliptic nodes),
26 with $\epsilon=-1$, 27 with $\epsilon=0,-2$,
or 29 with $\epsilon=\pm 2,0$,
by resolving them according to $c=\pm 1$.
\end{coro}
\begin{proof}
The proof is similar to that of Corollary
\ref{coro-2ptg0}. Note that if the number of
odd arcs is not more than 4 then we have the moves
of Proposition \ref{prop-moves}
also for the genus 1 case.
\end{proof}

\subsection{Trinodal spatial elliptic sextics
and proof of Theorems \ref{thm-g1t} and \ref{thm-g1r}}
Additional isotopies among the curves corresponding
to positive and negative resolution of nodal
diagrams from Figure \ref{g1-2pt} are obtained
with the help of bi- and trinodal spatial elliptic
sextics. As in the case of quadrinodal rational
sextics we may resolve the nodes of such curves
independently according to our choice of signs.
%Furthermore we may smooth a node in any quadrinodal
%rational sextic to obtain a trinodal elliptic sextic
%with the help of the following proposition.

\begin{lem}\label{34red}
Let $J\subset\p^3$ be a rational real quadrinodal sextic
and $p\in J$ be one of its real nodes (can be elliptic
or hyperbolic). %For every choice of real smoothing of $p$
We may choose to smooth $J$ at $p$
to a real elliptic sextic $I$ of a type I or to a type II
keeping three other nodes of $J$.
\end{lem}
\begin{proof}
We have seen that $J$ is given by 4 chords on a conic
$Q\subset\p^2$.
We can view the chords as lines in $\p^2$.
Perturbing the diagram if needed
we may assume that no three chords intersect
in a point. 
The plane $p^2$ can be linearly embedded to $\p^3$
so that the four chords are cut by the coordinate
planes. Then $J$ is the image of $Q$ under \eqref{cucr}.

Let the line $L_p\subset\p^2$ be the chord corresponding to $p$. Let $R\subset\p^2$ be the cubic curve obtained
by perturbation of $J\cup L_p$ with the help of the three
other chords (either to a type I or type II
real curve). The image of $R$ under \eqref{cucr}
is $I$.
%Then the four chords on $\p^2$ define
%four hyperplane sections of 
% can be viewed
%as a line $L\subset\p^2$.
\end{proof}

%\begin{proof}
\noindent{\em Proof of Theorems  \ref{thm-g1t} and \ref{thm-g1r}.}
%\end{proof}
Suppose that $K$ is of type II.
We use the same quadrinodal
curves as in Figure \ref{4dp} and apply to them
Lemma \ref{34red} to complete the classification
in this case.

\begin{figure}[h]
\includegraphics[width=100mm]{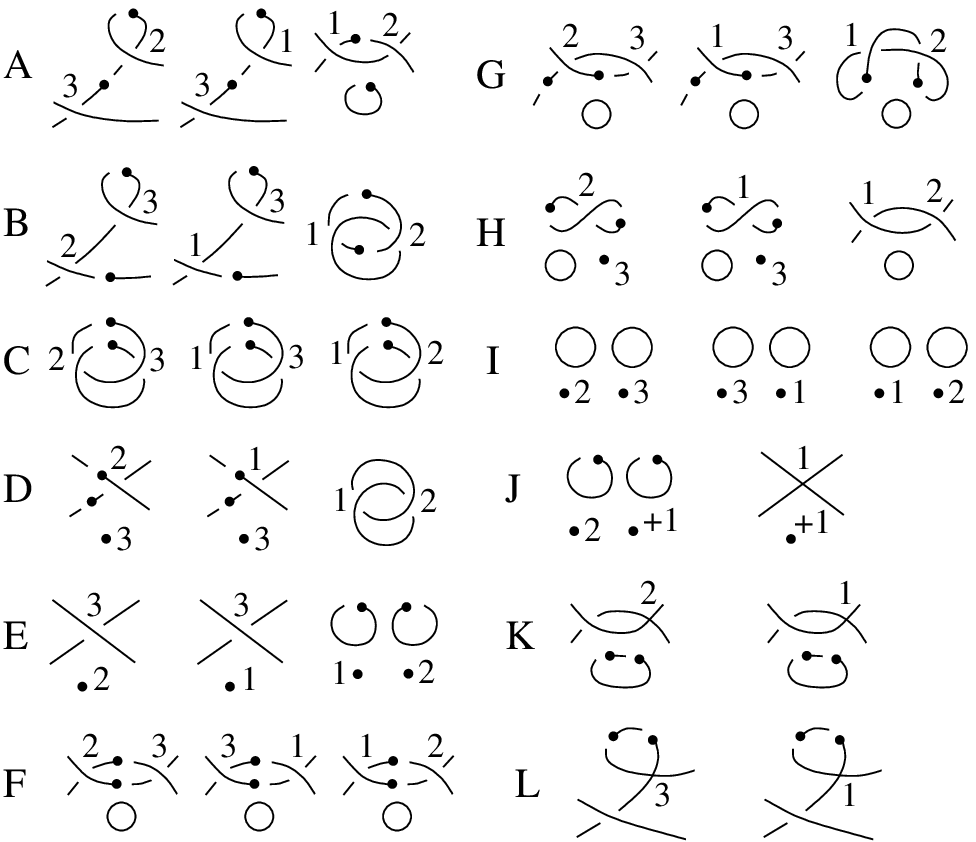}
\caption{
Twelve trinodal and binodal elliptic curves
\label{3dp}}
\end{figure}

\begin{figure}[h]
\includegraphics[width=125mm]{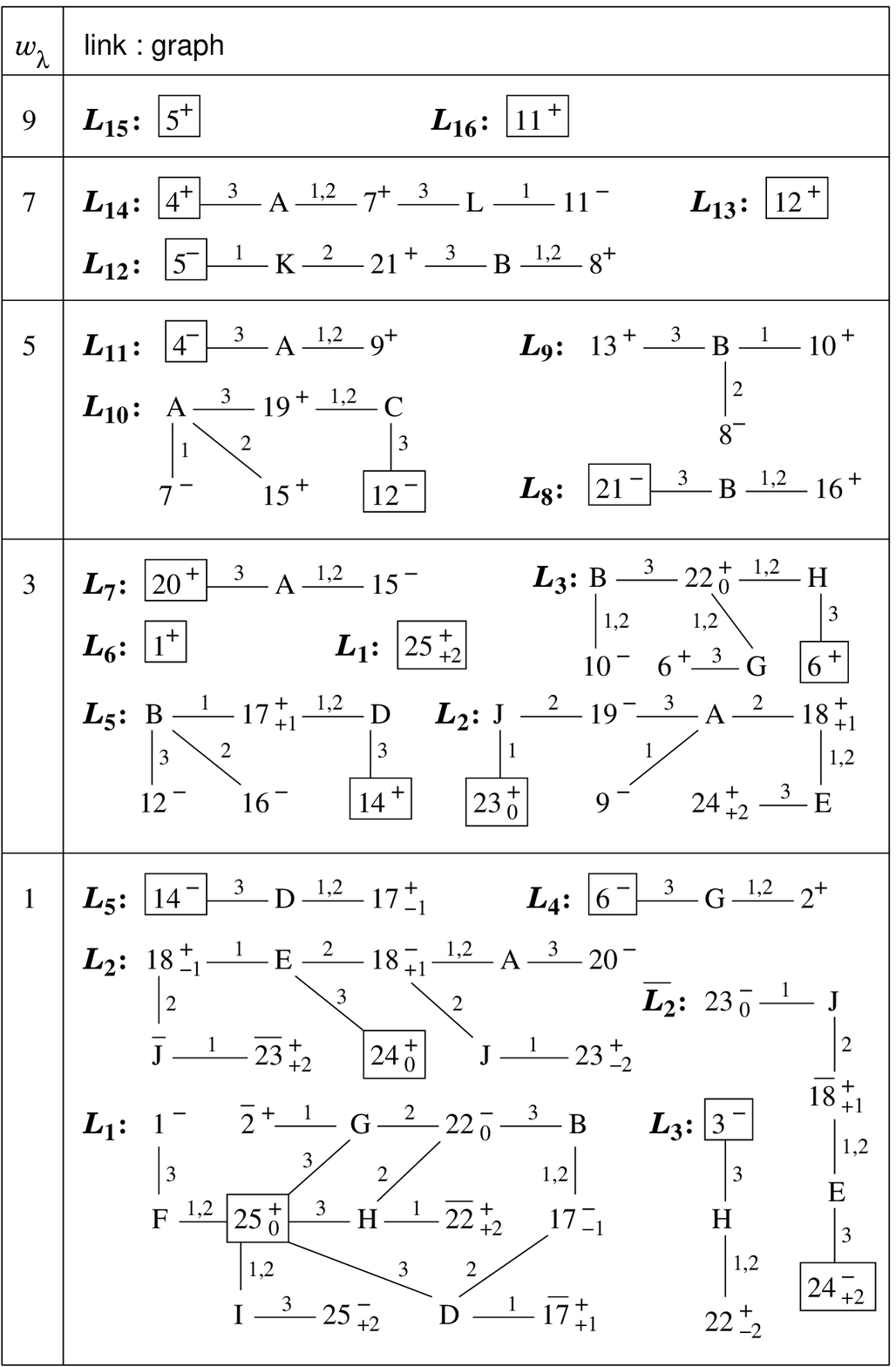}
\caption{
Graphs of rigid isotopy equivalence.
\label{graphs-d6g1}}
\end{figure}

If $K$ is of type I we use trivalent curves from
Figure \ref{3dp} to relate the curves from 
Corollary \ref{coro-2ptg1} with graphs depicted
at Figure \ref{graphs-d6g1}.
\ignore{
All the trinodal curves in Figure \ref{3dp}
except for J
can be obtained from quadrinodal curves
of Figure \ref{4dp} with the help of Lemma \ref{34red}.
It is easy to construct J as well
as the binodal curves K and L
explicitly.
}
All the nodal curves in Figure 15 except for J can be obtained from
rational nodal curves of Figures 9 and 10 with the help of Lemma 33
(see Table below).
It is easy to construct J explicitly.
\begin{center}
  \begin{tabular}{|c|c|c|c|c|c|c|c|c|c|c|c|}
    \hline
Rational multinodal curve & A & B & C & H &
J & H & I & J & K & L & L \\
\hline
Perturbed node & 1 & 2 & 2 & 2 & 3 & 1 & 4 &
1 & 1 & 3 & 1\\
\hline 
Elliptic multinodal curve & 
A & B & C & D & E & F & G & H & I & K & L\\
\hline
\end{tabular}
\end{center}


\qed

%% file: main.bbl
\begin{thebibliography}{99}
\bibitem{Bj}
{Bj{\"o}rklund, Johan},
   {\em Real algebraic knots of low degree},
   {J. Knot Theory Ramifications},
   {\bf 20:9} (2011),
   {1285--1309}. 
      
 \bibitem{DeMello}
 D'Mello, Shane,
 {\em Rigid isotopy classification of real degree-4 planar rational curves with only real nodes (An elementary approach)},
  arXiv:1307.7456.
   
   
 \bibitem{JuViro} {Drobotukhina, J.},
  {\em Classification of links in $RP^3$ with at most six crossings},
%in ``Topology of manifolds and varieties'' (ed. O.Ya. Viro),
Advances in Soviet Math., {\bf 18} (1994), 87--121.

\bibitem{Harnack}
   {Harnack, Axel},
   {\em \"Uber die Vieltheiligkeit der ebenen algebraischen Curven},
   {Math. Ann.},
   {\bf 10:2} (1876),
   {189--198}.
   
\bibitem{Hilbert} {Hilbert, David},
   {\em Mathematical problems},
   {Bull. Amer. Math. Soc. {\bf 8} (1902), 437--479}.
   
\bibitem{Rokhlin} Rokhlin, V. A.,
{\em Complex topological characteristics of real algebraic curves},
Uspehi Mat. Nauk, {\bf 33:5} (1978), 77-89.
 
\bibitem{Vi}
{Viro, Oleg},
   {\em Encomplexing the writhe},
{Amer. Math. Soc. Transl. Ser. 2},
      {\bf 202} (2001),
{241--256}.
  

\end{thebibliography}
